\documentclass[a4paper,final,oneside,notitlepage,11pt]{article}

\usepackage{amssymb}
\usepackage{amsmath}
\usepackage{amstext}
\usepackage[]{geometry}
\usepackage{mathrsfs}
\usepackage{euscript}
\usepackage[all]{xy}
\usepackage{xstring}
\usepackage{slashed}

\def\N {\mathbb{N}}
\def\Z {\mathbb{Z}}

\def\R {\mathbb{R}}

\def\id{\mathrm{id}}

\def\quand{\quad\text{ and }\quad}

\def\ev{\mathrm{ev}}
\def\lc{_{\mathrm{lc}}}

\def\hc#1{\mathrm{h}_{#1}}

\def\subset{\subseteq}

\def\nobr{~\hspace{-0.26em}}
\def\maps{\nobr:\nobr}
\def\df{\nobr := \nobr}

\let\Oldin\in\renewcommand{\in}{\nobr\Oldin\nobr}
\let\Oldtimes\times\renewcommand{\times}{\nobr\Oldtimes}
\let\Oldotimes\otimes\renewcommand{\otimes}{\nobr\Oldotimes}

\newlength{\widthtmp}
\def\length#1{\settowidth{\widthtmp}{#1}\the\widthtmp}

\def\lli#1{\,_{#1}\!}

\renewcommand{\varepsilon}{\epsilon}


\def\erf#1{(\ref{#1})}

\newlength{\myl}
\newcommand\sheaf[1]{\unitlength 0.1mm
  \settowidth{\myl}{$#1$}
  \addtolength{\myl}{-0.8mm}
  \begin{picture}(0,0)(0,0)
  \put(3,6){\text{\uline{\hspace{\myl}}}}
  \end{picture}#1\hspace{-0.15mm}}

\def\man{\mathcal{M}\!a\!n}
\def\diff{\mathcal{D}\!i\!f\!\!f}
\def\frech{\mathcal{F}\hspace{-0.12em}r\hspace{-0.1em}e\hspace{-0.05em}c\hspace{-0.05em}h}
\def\top{\mathcal{T}\!\!op}
\def\set{\mathcal{S}\hspace{-0.1em}et}

\def\hom{\mathcal{H}\!om}

\def\brackets#1{\IfStrEq{#1}{-}{}{(#1)}}
\def\subindex#1{\IfStrEq{#1}{-}{}{_{#1}}}
\def\buntech#1#2{\mathcal{B}\hspace{-0.01em}un_{\hspace{-0.1em}#1}^{#2}}
\def\bun#1#2{\buntech{#1}{}\brackets{#2}}
\def\buncon#1#2{\buntech{#1}{\nabla}\hspace{-0.05em}\brackets{#2}}
\def\bunconflat#1#2{\buntech#1{\nabla_{\!0}}\hspace{-0.05em}\brackets{#2}}

\newcommand{\alxydim}[2]{\begin{aligned}\xymatrix#1{#2}\end{aligned}}

\renewcommand{\to}{\nobr\!\xymatrix@R=0cm@C=1.4em{\ar[r] &}\nobr}
\renewcommand{\mapsto}{\!\xymatrix@R=0cm@C=1.4em{\ar@{|->}[r] &}\!}
\renewcommand{\Rightarrow}{\!\xymatrix@R=0cm@C=1.4em{\ar@{=>}[r] &}\!}
\renewcommand{\Leftarrow}{\!\xymatrix@R=0cm@C=1.4em{\ar@{<=}[r] &}\!}
\newcommand{\incl}{\!\xymatrix@R=0cm@C=1.4em{\ar@{^(->}[r] &}\!}

\renewcommand\Leftrightarrow{\!\xymatrix@R=0cm@C=1.4em{\ar@{<=>}[r] &}\!}

\usepackage{amsthm}
\usepackage{graphicx}
\usepackage[margin=2cm,font=normalsize,labelfont=bf]{caption}
\usepackage{verbatim}
\usepackage{enumerate}
\usepackage{fancyhdr}
\usepackage[]{geometry}
\usepackage[normalem]{ulem}
\usepackage{xstring}
\usepackage{needspace}

\makeatletter

\newcounter{denseversion}
\newcounter{authorcounter}
\newcounter{adresscounter}

\def\title#1{\gdef\@title{#1}}
\def\@title{}
\def\subtitle#1{\gdef\@subtitle{#1}}
\def\@subtitle{}

\def\authortagsused{0}
\def\adresstag#1{\if!#1!\else$^{\;#1\;}$\fi}

\renewcommand{\author}[2][]{
  \stepcounter{authorcounter}
  \if!#1!\else\gdef\authortagsused{1}\fi
  \ifnum\value{authorcounter}=1
    \def\@authorstringa{#2\adresstag{#1}}
    \def\@authorstringb{#2}
    \def\@authorstringc{#2\adresstag{#1}}
  \else
    \g@addto@macro\@authorstringa{\ and #2\adresstag{#1}}
    \g@addto@macro\@authorstringb{\ and #2}
    \g@addto@macro\@authorstringc{\\#2\adresstag{#1}}
  \fi}
\def\@author{\ifnum\value{denseversion}=0\@authorstringa\else\@authorstringb\fi}

\def\@adressstringa{}
\def\@adressstringb{}
\newcommand{\adress}[2][]{
  \stepcounter{adresscounter}
  \ifnum\value{adresscounter}=1
    \g@addto@macro\@adressstringa{\ifnum\authortagsused=0\def\br{\\}\else\def\br{, }\fi\adresstag{#1}#2}
    \g@addto@macro\@adressstringb{\def\br{\\}\adresstag{#1}\parbox[t]{14cm}{#2}}
  \else
    \g@addto@macro\@adressstringa{\\[\bigskipamount]\adresstag{#1}#2}
    \g@addto@macro\@adressstringb{\\[\medskipamount]\adresstag{#1}\parbox[t]{14cm}{#2}}
  \fi}
\def\@adress{\ifnum\value{denseversion}=0\@adressstringa\else\@adressstringb\fi}

\def\preprint#1{\gdef\@preprint{#1}}
\def\@preprint{}
\def\keywords#1{\gdef\@keywords{#1}}
\def\@keywords{}
\def\msc#1{\gdef\@msc{#1}}
\def\@msc{}
\def\email#1{
   \gdef\@email{#1}
   \g@addto@macro\@authorstringc{ {\it (#1)}}}
\def\@email{}
\def\dedication#1{\gdef\@dedication{#1}}
\def\@dedication{}

\def\mybaselinestretch#1{\gdef\@mybaselinestretch{#1}}
\def\@mybaselinestretch{}

\def\refname{References}

\mybaselinestretch{1.2}
\renewcommand{\baselinestretch}{\@mybaselinestretch}

\def\denseversion{
  \setcounter{denseversion}{1}
  \newgeometry{left=3cm,right=3cm,top=3cm}
  \mybaselinestretch{1.1}
  \renewcommand{\baselinestretch}{\@mybaselinestretch}
  \normalfont
  \fancyfoot[C]{\itshape{\hspace{2.5cm}--$\,\,$\thepage$\,\,$--}}}

\newlength{\myparskip}
\newlength{\myproofparskip}

\renewcommand{\emph}[1]{\def\reserved@a{it}\ifx\f@shape\reserved@a\uline{#1}\else\textit{#1}\fi}

\newcommand{\mytableofcontents}{
   \ifnum\value{denseversion}=0
     \tableofcontents
   \else
     \renewcommand{\baselinestretch}{0.8}
     \normalfont
     \tableofcontents
     \renewcommand{\baselinestretch}{\@mybaselinestretch}
     \normalfont
   \fi}

\newlength{\zeilenlaenge}
\def\putindent#1{
  \settowidth{\zeilenlaenge}{#1}
  \ifnum\zeilenlaenge>\textwidth
    #1
  \else
    \noindent #1
  \fi
}

\def\href#1#2{#2}

\def\kohyp{
  \usepackage{hyperref}
  \hypersetup{
    linktocpage = true,
    pdftitle = {\@title},
    pdfauthor = {\@author},
    pdfkeywords = {\@keywords},    
    bookmarksopen = true,
    bookmarksopenlevel = 1
  }}  
\def\showkeywords{\begin{flushleft}\footnotesize\textbf{Keywords}: \@keywords\end{flushleft}}
\def\showmsc{\begin{flushleft}\footnotesize\textbf{MSC 2010}: \@msc\end{flushleft}}

\newcounter{mythm}[subsection]
\newcounter{mainthm}

\def\setsecnumdepth#1{
  \setcounter{secnumdepth}{#1}
  \setcounter{mythm}{0}
  \ifnum \c@secnumdepth >0
    \ifnum \c@secnumdepth >1
      \def\themythm{\thesubsection.\arabic{mythm}}
      \numberwithin{equation}{subsection}
      \renewcommand\theequation{\thesubsection.\arabic{equation}}
    \else
      \def\themythm{\thesection.\arabic{mythm}}
      \numberwithin{equation}{section}
      \renewcommand\theequation{\thesection.\arabic{equation}}
    \fi
  \else
    \def\themythm{\arabic{mythm}}
  \fi}

\newenvironment{mythmenv}{\strut\ \setlength{\parskip}{\myproofparskip}}{\setlength{\parskip}{\myparskip}}

\newlength{\mythmskip}
\newlength{\mythmtopskip}
\newtheoremstyle{mythmstylea}{\mythmtopskip}{\mythmskip}{\it}{}{\bf}{.}{0em}{}
\newtheoremstyle{mythmstyleb}{\mythmtopskip}{\mythmskip}{}{}{\bf}{.}{0em}{}

\theoremstyle{mythmstylea}
\newtheorem{mytheorem}[mythm]{Theorem}
\newtheorem{mydefinition}[mythm]{Definition}
\newtheorem{mycorollary}[mythm]{Corollary}
\newtheorem{myproposition}[mythm]{Proposition}
\newtheorem{mylemma}[mythm]{Lemma}

\newtheorem{mymaintheorem}[mainthm]{Theorem}
\newtheorem{mymaincorollary}[mainthm]{Corollary}
\newtheorem{mymainproposition}[mainthm]{Proposition}
\newtheorem{mymaindefinition}[mainthm]{Definition}

\theoremstyle{mythmstyleb}
\newtheorem{myremark}[mythm]{Remark}
\newtheorem{myexample}[mythm]{Example}
\newtheorem{myexercise}[mythm]{Exercise}

\newenvironment{theorem}[1][]{\begin{mytheorem}[#1]\begin{mythmenv}}{\end{mythmenv}\end{mytheorem}}
\newenvironment{definition}[1][]{\begin{mydefinition}[#1]\begin{mythmenv}}{\end{mythmenv}\end{mydefinition}}

\newenvironment{proposition}[1][]{\begin{myproposition}[#1]\begin{mythmenv}}{\end{mythmenv}\end{myproposition}}
\newenvironment{lemma}[1][]{\begin{mylemma}[#1]\begin{mythmenv}}{\end{mythmenv}\end{mylemma}}
\newenvironment{remark}[1][]{\begin{myremark}[#1]\begin{mythmenv}}{\end{mythmenv}\end{myremark}}
\newenvironment{example}[1][]{\begin{myexample}[#1]\begin{mythmenv}}{\end{mythmenv}\end{myexample}}

\newenvironment{maintheorem}[1]{\def\themainthm{#1}\begin{mymaintheorem}\begin{mythmenv}}{\end{mythmenv}\end{mymaintheorem}}

\newenvironment{maincorollary}[1]{\def\themainthm{#1}\begin{mymaincorollary}\begin{mythmenv}}{\end{mythmenv}\end{mymaincorollary}}

\renewenvironment{proof}[1][Proof]{\noindent #1. \begin{mythmenv}}{\hfill$\square$\end{mythmenv}\medskip}

\def\mytitle{}
\def\zmptitle{
  \begin{tabular}{cc}
    \begin{minipage}[c]{0.4\textwidth}
      \begin{flushleft}
        \includegraphics[width=110pt]{../../tex/zmp}
      \end{flushleft}  
    \end{minipage}&
    \begin{minipage}[c]{0.55\textwidth}
      \begin{flushright}
      {\small\sf\@preprint}
      \end{flushright}
    \end{minipage}
  \end{tabular}
  \vskip 2cm}

\def\maketitle{
  \setlength{\parskip}{\myparskip}  
  \newpage
  \noindent
  \mytitle
  \begin{center}
    \LARGE\@title\\
    \if!\@subtitle!\else \smallskip\LARGE\@subtitle\\\fi
    \bigskip
    \if!\@author!\else\bigskip\large\@author\\\fi
    \ifnum\value{denseversion}=0
      \if!\@adress!\else     \bigskip\normalsize\@adress\\\fi
      \if!\@email!\else\ifnum\value{authorcounter}=1\bigskip\normalsize\textit{\@email}\\\else\fi\fi
    \else
    \fi
    \if!\@dedication!\else \bigskip\normalsize{\@dedication}\\\fi
  \end{center}
  \ifnum\value{denseversion}=0\vskip 1.5cm\else\vskip0.5cm\fi
  \thispagestyle{empty}}

\def\kobiburl#1{
   \IfBeginWith
     {#1}
     {http://arxiv.org/abs/}
     {\kobibarxiv{#1}}
     {\kobiblink{#1}}}
\def\kobibarxiv#1{\href{#1}{\texttt{[arxiv:\StrGobbleLeft{#1}{21}]}}}
\def\kobiblink#1{Available as: \href{#1}{\texttt{\StrSubstitute{#1}{_}{\underline{\;\;}}}}}

\def\kobib#1{
  \begin{raggedright}
  \ifnum\value{denseversion}=0\else\small\fi

  \end{raggedright}
  \ifnum\value{denseversion}=0\else
      \noindent
      \if!\@authorstringc!\else
        \ifnum\authortagsused=0\ifnum\value{authorcounter}>1\normalsize\@authorstringc\\[\medskipamount]\else\fi\else\normalsize\@authorstringc\\[\medskipamount]\fi
      \fi
      \if!\@adress!\else\normalsize\@adress\\{}\fi
      \ifnum\authortagsused=0
        \ifnum\value{authorcounter}=1
          \if!\@email!\else\linebreak\normalsize\textit{\@email}\\{}\fi
        \else
        \fi
      \else
      \fi
  \fi
  }

\newenvironment{commentfigure}{\begin{comment}}{\end{comment}}
\newenvironment{sidewayscommentfigure}{\begin{minipage}}{\end{minipage}}

\def\showcomments{ -- Comments suppressed}

\newif\if@fewtab\@fewtabtrue{
  \count255=\time\divide\count255 by 60
  \xdef\hourmin{\number\count255}
  \multiply\count255 by-60\advance\count255 by\time
  \xdef\hourmin{\hourmin:\ifnum\count255<10 0\fi\the\count255}}
\def\ps@draft{
  \let\@mkboth\@gobbletwo
  \def\@oddfoot{
    \hbox to 7 cm{\tiny \versionno\hfil}
    \hskip -7cm\hfil\rm\thepage\hfil{\tiny\draftdate}}
  \def\@oddhead{}
  \def\@evenhead{}
  \let\@evenfoot\@oddfoot}
\def\draftdate{\number\month/\number\day/\number\year\ \ \ \hourmin }
\newcommand\version[1]{
  \typeout{}\typeout{#1}\typeout{}
  \vskip-1.7cm \centerline{\fbox{{\normalsize\tt DRAFT -- #1 -- 
  \draftdate\showcomments}}} \vskip0.92cm}
\def\draft#1{
  \def\versionno{#1}
  \pagestyle{draft}\thispagestyle{draft}
  \gdef\@ntitle{\version\versionno \@title}
  \global\def\draftcontrol{1}}
\global\def\draftcontrol{0}

\makeatother

\usepackage[latin1]{inputenc}
\usepackage[english]{babel}

\def\quot#1{``#1''}

\def\exd#1{{#1^{\vee}}}

\def\pcl#1{P\!#1_{c\hspace{0.01em}l}}

\def\ptx#1#2{\mathcal{P}_{\!#2}#1}
\def\ev{\mathrm{ev}}
\def\pt#1{\mathcal{P}#1}
\def\lt#1{\mathcal{L}#1}

\def\hc#1{\mathrm{h}_{#1}}
\def\pcomp{\nobr\star\nobr}

\def\prev#1{\overline{#1}}

\def\fun#1#2{\mathcal{F}\!un(#1,#2)}

\def\un{\mathscr{R}}
\def\uncon{\mathscr{R}^{\nabla}}
\def\tr{\mathscr{T}}
\def\trcon{\mathscr{T}^{\nabla}}

\def\fus#1#2{\mathcal{F}\!us(#2,#1)}
\def\fuslc#1#2{\mathcal{F}\!us_{\mathrm{lc}}(#2,#1)}
\def\fushom#1#2{h\!\fus {#1} {#2}}

\def\diffbun#1#2{\diff\hspace{-0.1em}\bun{#1}{#2}}
\def\diffbuncon#1#2{\diff\hspace{-0.1em}\buncon{#1}{#2}}
\def\diffbunconflat#1#2{\diff\hspace{-0.1em}\bunconflat#1{#2}}

\def\ptr#1{\tau_{#1}}
\def\ptrcon#1#2{\tau_{#1}^{\,#2}}

\title{Transgression to Loop Spaces and its Inverse, I}
\subtitle{Diffeological Bundles and Fusion Maps}
\author{Konrad Waldorf}
\adress{Department of Mathematics\\University of California, Berkeley\\970 Evans Hall \#3840\\Berkeley, CA 94720, USA}
\email{waldorf@math.berkeley.edu}
\keywords{bundle, connection, diffeological space, holonomy, loop space, thin homotopy, transgression}
\msc{Primary 53C29, Secondary 58B25}

\kohyp
\usepackage{amstext}
\usepackage[normalem]{ulem}

\begin{document}

%\draft{Version 2.0, 31.3.2011}
%\comments

\maketitle 

\begin{abstract}
We prove that isomorphism classes of  principal bundles  over a diffeological space are in bijection to certain maps on its  free loop space, both in a setup with and without connections on the bundles.  The maps on the loop space are smooth and satisfy a \quot{fusion} property with respect to triples of paths. Our bijections are established by explicit group isomorphisms: transgression and regression. Restricted to smooth, finite-dimensional manifolds, our results extend previous work of J. W. Barrett. 
\showkeywords
\showmsc
\end{abstract}

\tableofcontents

\section{Introduction and Results}

\label{sec:intro}

We study a relationship between geometry on a space and geometry on its  loop space. We are concerned with a fairly general class of spaces: diffeological spaces, one version of the \quot{convenient calculus} \cite{souriau1,kriegl1}. Most prominently, the category of diffeological spaces contains the categories of smooth manifolds and Fréchet manifolds as full subcategories, and a lot of familiar geometry generalizes almost automatically from these subcategories  to diffeological spaces.

The geometry we study on  a diffeological space $X$ consists of principal bundles with connection. The structure group of these bundles is an ordinary, abelian Lie group $A$, which is allowed to be discrete and non-compact. The loop space $\lt X$ that is relevant here consists of so-called thin homotopy classes of smooth maps $\tau\maps S^1 \to X$, and is thus better called the \quot{thin loop space} of $X$. Due to the convenient properties of diffeological spaces, $\lt X$ is again an honest diffeological space. On  $\lt X$, we characterize a class of \quot{fusion maps} $f\maps  \lt X \to A$, following an idea of Stolz and Teichner \cite{stolz3}. A fusion map $f$ is smooth, and whenever  $\gamma_1$, $\gamma_2$ and $\gamma_3$ are paths in $X$ with a common initial point and a common end point, it satisfies
\begin{equation*}
f(\prev {\gamma_2} \pcomp \gamma_1) \cdot f(\prev{\gamma_3} \pcomp \gamma_2) = f(\prev{\gamma_3} \pcomp \gamma_1)\text{,}
\end{equation*}
where $\pcomp$ denotes the composition of paths and $\prev\gamma$ denotes the reversed path. 
Fusion maps form a group under point-wise multiplication, which we denote by  $\fus A{\lt X}$. A detailed discussion of fusion maps is the content of Section \ref{sec:pathsloops}. Our first result is
\begin{maintheorem}{A}
\label{th:main}
Let $X$ be a connected diffeological space and let $A$ be an abelian Lie group. There is an isomorphism   
\begin{equation*}
\hc 0 \diffbuncon A X \cong \fus A {\lt X}
\end{equation*}
between the group of isomorphism classes of diffeological principal $A$-bundles over $X$ with connection and the group of fusion maps on the thin loop space $\lt X$.  
\end{maintheorem}

The bijection of Theorem \ref{th:main} is established by group isomorphisms called \quot{transgression} and \quot{regression}. Transgression basically takes the holonomy of the given connection, and is discussed in Section \ref{sec:transgression}.  Its inverse, regression, reconstructs a principal bundle over $X$ with connection from a given fusion map $f\maps \lt X \to A$, and is the content of Section \ref{sec:regression}. The proof of Theorem \ref{th:main} is given in  Section \ref{sec:proofs}. 

The groups on both  sides of the isomorphism of Theorem \ref{th:main} have important subgroups: the one on the left hand side is composed of bundles with \emph{flat} connections, and the one on the right hand side is composed of \emph{locally constant} fusion maps. In Propositions \ref{prop:flatreg} and \ref{prop:holder} we shall show the following:

\begin{maincorollary}{A}
\label{co:flat}
The isomorphism of Theorem \ref{th:main} restricts to an isomorphism
\begin{equation*}
\hc 0 \diffbunconflat A X \cong \fuslc A {\lt X}
\end{equation*}
between the group of isomorphism classes of \uline{flat} diffeological principal $A$-bundles over $X$ and the group of \uline{locally constant} fusion maps on  $\lt X$.
\end{maincorollary}

If $X$ is a smooth manifold, a \emph{diffeological} principal bundle over $X$ is the same as a \emph{smooth} principal bundle; similarly, connections on diffeological principal bundles become ordinary, smooth connections. 
Under these identifications, Theorem \ref{th:main} and Corollary \ref{co:flat} are statements about the  geometry of ordinary,  smooth manifolds. 

Our second result is an analog of Theorem \ref{th:main} in a setup for bundles \emph{without} connection. It relies on the existence of connections on principal bundles and is thus only valid over smooth manifolds, and not over general diffeological spaces. The relevant structure on the thin loop space is now  a group of \emph{equivalence classes} of fusion maps, which we denote by $h\fus A{\lt M}$. Here two fusion maps are identified if they are connected by a path through the  space of fusion maps. 

\begin{maintheorem}{B}
\label{th:man}
Let $M$ be a connected smooth manifold. There is an isomorphism
\begin{equation*}
\hc 0 \bun AM \cong \fushom A{\lt M}
\end{equation*}
between the group of isomorphism classes of smooth principal $A$-bundles over $M$ and the group of equivalence classes of fusion maps on $\lt M$.
\end{maintheorem}

Our last result is that under the isomorphisms of Theorems \ref{th:main} and \ref{th:man}, the operation of forgetting connections corresponds precisely to the projection of a fusion map to its equivalence class:
\begin{maintheorem}{C}
\label{th:comp}
The bijections of Theorems \ref{th:main} and \ref{th:man} fit into a commutative diagram
\begin{equation*}
\alxydim{@=1.3cm}{\hc 0 \buncon A M \ar@{<->}[r]^-{\cong} \ar[d] & \fus A {\lt M} \ar[d] \\ \hc 0 \bun A M \ar@{<->}[r]_-{\cong} & \fushom A {\lt M}\text{,}}
\end{equation*}
whose vertical arrows are, respectively, forgetting the connection and projecting to equivalence classes. 
\end{maintheorem}

Theorems \ref{th:man} and \ref{th:comp} are proved simultaneously. In Sections \ref{sec:regression} and \ref{sec:transgression} we define group homomorphisms $\tr$ and $\un$ that constitute the isomorphism of Theorem \ref{th:man}  \emph{such that} Theorem \ref{th:comp} is true. In other words, $\tr$ and $\un$ are covered, respectively, by transgression and regression, which are inverses of each other according to Theorem \ref{th:main}. Since the vertical arrows in the commutative diagram of Theorem \ref{th:comp} are \emph{surjective}, it follows that  $\tr$ and $\un$ are also inverses to each other.

Studying the relationship between principal bundles with connection over a smooth manifold $M$ and their holonomy  has a long history in differential geometry, including work of Milnor, Lashof and Teleman. A nice overview and references  can be found in Barrett's seminal paper \cite{barret1}. In that paper, Barrett introduced the notion of thin homotopic loops that we use here, and proved a version of our Theorem \ref{th:main} for smooth manifolds.  We remark that his formulation of the loop space structure is slightly different from our fusion maps, and requires  fixing a base point in $M$.

The results of the present article extend Barrett's result in two ways. The first  is that Theorem \ref{th:main} extends Barrett's bijection to a larger class of spaces -- diffeological spaces. For example, Theorem \ref{th:main} holds for principal bundles over the thin loop space itself. The second  is  Theorem \ref{th:man}, which extends Barrett's bijection to principal bundles \emph{without} connection.

Our  results  are made possible due to new tools that we develop. The main innovation is  a comprehensive theory of principal bundles with connection over diffeological spaces; this is worked out in Section \ref{sec:diffbun}. We follow the slogan: even if one is only interested in smooth manifolds it is  helpful to use diffeological spaces.  More specifically, we use new  descent-theoretical aspects:  we introduce a Grothendieck topology on the category of diffeological spaces and prove that principal bundles with and without connections form sheaves of groupoids (Theorems \ref{th:diffbunsheaf} and \ref{th:diffbunconsheaf}). Finally, we show that \emph{diffeological} principal bundles with and without connection reduce over smooth manifolds consistently to  \emph{smooth} principal bundles (Theorems \ref{th:equivalencesmooth} and \ref{th:conequivalencesmooth}).

In the main text of this paper we assume that the reader is familiar with the basics of diffeological spaces. In order to make the paper accessible for others, we have included Appendix \ref{sec:diff} with a brief review about diffeological spaces, emphasizing  differential forms, path spaces and sheaf theory. 

\paragraph{Acknowledgements.} I gratefully acknowledge  a Feodor-Lynen scholarship, granted by the Alex\-an\-der von Hum\-boldt Foundation. Further, I thank the Max-Planck-Institut für Mathematik in Bonn for kind hospitality and support.  I thank Thomas Nikolaus, Martin Olbermann, Arturo Prat-Waldron, Urs Schreiber, Andrew Stacey, and Peter Teichner for many exciting discussions.

\setsecnumdepth{2}

\section{Paths and Loops in Diffeological Spaces}

\label{sec:pathsloops}

In this section we define fusion maps on thin loop spaces. They appear on the right hand side of the bijections of Theorems \ref{th:main} and \ref{th:man}. 
\subsection{Thin Homotopy Equivalence}

\label{sec:thinhom}

Let $X$ be a diffeological space. A \emph{path} in $X$ is a smooth map $\gamma\maps [0,1] \to X$ that is locally constant in a neighborhood of $\left \lbrace 0,1 \right \rbrace$. The latter condition is also known under the name  \quot{sitting instants}. Our main reference for this section is \cite{iglesias1}, where paths are called \quot{stationary paths}. The set of paths is denoted $PX$; it is itself a diffeological space as a subset of the diffeological space of smooth maps from $[0,1]$ to $X$, which we denote by  $D^\infty([0,1],X)$.  

A diffeological space $X$ is called \emph{connected} if the endpoint evaluation 
\begin{equation*}
\ev\maps  PX \to X \times X\maps  \gamma \mapsto (\gamma(0),\gamma(1))
\end{equation*}
is surjective. One can show that $\ev$ is then even a subduction \cite[V.6]{iglesias1}, the diffeological analog of a smooth map that admits  smooth local sections (see Definition \ref{def:subduction} and Lemma \ref{lem:smoothsubduction}). 
In the following we assume that $X$ is connected.

Because of the sitting instants, two paths $\gamma_1,\gamma_2$ with $\gamma_1(1)=\gamma_2(0)$ can be composed to a third path $\gamma_2 \pcomp \gamma_1$ which is defined in the usual way \cite[V.2, V.4]{iglesias1}. For $x\in X$, we  denote by $\id_x$ the constant path at $x$, and for a path $\gamma$ we denote by $\prev\gamma$ the reversed path.
\begin{comment}
\begin{lemma}[{{\cite[V.3, V.4]{iglesias1}}}]
\label{lem:compsmooth}
Composition and reversion are smooth maps
\begin{equation*}
\pcomp\maps  PX \lli{\ev_0}\times_{\ev_1} PX \to PX
\quand
\prev{(..)}\maps  PX \to PX\text{,}
\end{equation*}
where $\ev_i \df  \mathrm{pr}_i \circ \ev\maps  PX \to X$ for $i=1,2$.
\end{lemma}
\end{comment}
A smooth map $f\maps X \to Y$ between diffeological spaces induces a smooth map 
\begin{equation*}
Pf\maps PX \to PY\maps  \gamma \mapsto f \circ \gamma
\end{equation*}
that takes composition and reversal of paths in $X$ to those of paths in $Y$ \cite[I.59]{iglesias1}. 

We want to force composition to be associative, the constant paths to be identities, and path reversal to provide inverses for this composition. One solution would be to identify homotopic paths. However, as  is known  in the case of manifolds,  a lot of the geometry will be lost: one better restricts to homotopies of rank  one, so-called \emph{thin} homotopies \cite{barret1,schreiber3}. In the following we generalize the concept of  thin homotopies to diffeological spaces. 

\begin{definition}
\label{def:thin}
Let $k\in \N$.
A smooth map $f\maps X \to Y$ between diffeological spaces has \emph{rank $k$} if for every plot $c\maps  U \to X$ and every point $x\in U$ there exists an open neighborhood $U_x \subset U$, a plot $d\maps V \to Y$ and a smooth map $g\maps  U_x \to V$ such that the diagram
\begin{equation*}
\alxydim{@=1.2cm}{U_{x} \ar[r]^{g} \ar[d]_{c} & V \ar[d]^{d} \\ X \ar[r]_{f} & Y}
\end{equation*}
is commutative, and the rank of the differential of $g$ is at most $k$. 
\end{definition}

With a view to Theorem \ref{th:man} it is important to see that notions we introduce for diffeological spaces reduce to the corresponding existent notions for smooth manifolds. 

\begin{lemma}
\label{lem:thinman}
For $M$ and $N$ smooth manifolds, a smooth map $f\maps M \to N$ has rank $k$ in the sense of Definition \ref{def:thin} if and only if its differential is at  most of rank $k$.
\end{lemma}

The proof is elementary. One can also show that a smooth map $f\maps X \to Y$ between diffeological spaces has rank $k$ if and only if Laubinger's tangential map $f_{*}$  has at most rank $k$ \cite{laubinger1}.
The following lemma summarizes obvious statements. 

\begin{lemma}
\label{lem:thinprop}
The rank of smooth maps satisfies the following rules:
\begin{enumerate}
\item[(a)]
If $k>l$, every rank $l$ map also has  rank $k$.

\item[(b)]
Every constant map has rank zero. 

\item[(c)]
If $f\maps  X \to Y$ has rank $k$ and descends along a subduction $p\maps  X \to Z$, then the quotient map $f'\maps Z \to Y$ also has rank $k$.

\item[(d)]
If $f\maps X \to Y$ has rank $k$, and $g\maps  W \to X$ and $h\maps  Y \to Z$ are smooth maps, then $h \circ f \circ g\maps  W \to Z$ also has rank $k$. \end{enumerate}
\end{lemma}

Using the notion of the rank of a smooth map, we define a suitable equivalence relation on the space $PX$ of paths in $X$. 

\begin{definition}
\label{def:thinhomotopy}
Let $\gamma_1$ and $\gamma_2$ be paths in $X$ with a common initial point $x$ and a common end point $y$.
A \emph{homotopy} between $\gamma_1$ and $\gamma_2$ is a path $h\in PPX$ with 
\begin{equation*}
\ev(h) = (\gamma_1,\gamma_2)
\quand
\ev(h(s))=(x,y)
\end{equation*}
for all $s\in[0,1]$.  A homotopy $h$ is called \emph{thin}, if the adjoint map
\begin{equation*}
\exd h\maps  [0,1]^2 \to X \maps  (s,t) \mapsto h(s)(t)
\end{equation*}
has rank one.
\end{definition}

\begin{comment}
\begin{remark}
We remark that in the category of diffeological spaces the adjunction \begin{equation*}
\exd {(..)} \maps  D^\infty(X,D^\infty(Y,Z)) \cong D^\infty(X \times Y,Z)
\end{equation*}
we have used in Definition \ref{def:thinhomotopy} is  a diffeomorphism
for arbitrary diffeological spaces $X$, $Y$ and $Z$ \cite[I.60]{iglesias1}. \end{remark}
\end{comment}

An important example of a thin homotopy is an orientation-preserving reparameterization;  or even any smooth map $\eta\maps  [0,1] \to [0,1]$ with $\eta(0)=0$ and $\eta(1)=1$. A  homotopy between $\gamma \circ \eta$  and $\gamma$ can be obtained from  a smooth  homotopy between $\eta$ and the identity $\id_{[0,1]}$. It is thin due to Lemmata \ref{lem:thinman} and \ref{lem:thinprop} (d).

\begin{lemma}
\label{lem:thinhomeq}
Being thin homotopic is an equivalence relation on the diffeological space $PX$ of paths in $X$. 
\end{lemma}

\begin{proof}
First of all, the identity $\id_{\gamma}\in PPX$ is thin, since  $\exd{\id_{\gamma}}$ factors through $[0,1]$, and thus has rank one by Lemmata  \ref{lem:thinman} and \ref{lem:thinprop} (d). If $h\in PPX$ is a thin homotopy, $\prev h$ is also thin since $\exd{(\prev h)}$ factors through $\exd h$. Finally, we have to show that if $h_1 \in PPX$ and $h_2 \in PPX$ are composable thin homotopies, the composition $h_2 \pcomp h_1$ is again thin. 
For this purpose we recall that the path composition  $h_2 \pcomp h_1\maps  [0,1] \to PX$ is defined using the subduction $U\df  U_1 \sqcup U_{\varepsilon} \sqcup U_2 \to [0,1]$,  where 
\begin{equation*}
\textstyle
U_1 \df  [0,\frac{1}{2})
\quad\text{, }\quad
U_{\varepsilon} \df  (\frac{1}{2}-\frac{1}{2}\varepsilon,\frac{1}{2}+\frac{1}{2}\varepsilon)
\quand
U_2 \df  (\frac{1}{2},1]\text{,}
\end{equation*}
and $\varepsilon$ is chosen such that $h_1(t)=h_1(1)$ for all $t>1-\varepsilon$ and $h_2(t)=h_2(0)$ for all $t<\varepsilon$. We define $\tilde h\maps  U \to PX$ by $h_1(2t)$ over $U_1$, $h_1(1)=h_2(0)$ over $U_{\varepsilon}$ and $h_2(2t-1)$ over $U_2$. Obviously, $\tilde h$ descends to $[0,1]$; this defines $h_2 \pcomp h_1$ (see Lemma \ref{lem:quotmap}). Since $h_1$ and $h_2$ are thin, it follows that $\exd {\tilde h}\maps  U \times [0,1] \to X$ has rank one. Furthermore, it descends to $\exd{(h_2 \pcomp h_1)}$. Thus, $\exd{(h_2 \pcomp h_1)}$ has rank one by Lemma \ref{lem:thinprop} (c).
\end{proof}

Denoting by $\sim$ the equivalence relation of being thin homotopic,
the diffeological space of thin homotopy classes of paths is  
\begin{equation*}
\pt X \df  PX / \sim\text{.}
\end{equation*}
It remains to show that $\pt X$ has the desired structure:

\begin{proposition}
\label{prop:compdesc}
Path composition and reversal  
  descend to smooth maps
\begin{equation*}
\pcomp\maps  \pt X \times_X \pt X \to \pt X
\quand
\prev{(..)} \maps  \pt X \to \pt X\text{.}
\end{equation*}
The composition is associative and the constant paths are identities. Furthermore,  reversing paths provides inverses for the composition.
\end{proposition}

\begin{proof}
Composition and reversal of paths in $PX$ are smooth maps \cite[V.3, V.4]{iglesias1}. By Lemma \ref{lem:quotmap}, a smooth map that descends along a subduction descends to a \emph{smooth} map. So we only have to show that composition and reversal are well-defined under thin homotopies.

In order to do so, we introduce a \quot{pointwise} composition and reversal for paths in path spaces. If $h\in PPX$  is such a path, $r(h)\in PPX$ is defined by $r(h)(s) \df  \prev{h(s)}$. Since $\exd{r(h)}$ factors through $\exd h$, $r(h)$ is  thin whenever $h$ is thin. Thus, if  $h$ is a thin homotopy between $\gamma_1$ and $\gamma_2$, then $r(h)$ is a thin homotopy between $\prev{\gamma_1}$ and $\prev{\gamma_2}$.

The composition goes similarly; here one constructs from two paths $h_1,h_2\in PPX$ such that $h_1 \times h_2$ is a path in $PX \times_X PX$ a new path $c(h_1,h_2)$ by $c(h_1,h_2)(s) \df  h_2(s)\pcomp h_1(s)$. This is  thin  by the same reasoning as in the proof of Lemma \ref{lem:thinhomeq}. Summarizing, if $h_1$ is a thin homotopy between $\gamma_1$ and $\gamma_1'$, and $h_2$ is a thin homotopy between $\gamma_2$ and $\gamma_2'$, where $\gamma_1$ and $\gamma_2$ are composable, $c(h_1,h_2)$ is well-defined since
\begin{equation*}
h_1(s)(1) = \gamma_1(1) = \gamma_2(0) = h_2(s)(0)
\end{equation*}
for all $s\in[0,1]$,
and thus is a thin homotopy between $\gamma_2 \pcomp \gamma_1$ and $\gamma_2' \pcomp \gamma_1'$.

Composition is associative: one finds 
for three composable paths $\gamma_1$, $\gamma_2$ and $\gamma_3$,  a thin homotopy $(\gamma_3 \pcomp \gamma_2) \pcomp \gamma_1 \sim \gamma_3 \pcomp (\gamma_2 \pcomp \gamma_1)$
constructed from the evident reparameterization. In the same way one finds thin homotopies $\gamma \pcomp \id_x \sim \gamma \sim \id_y \pcomp \gamma$
for a path $\gamma$ with $\ev(\gamma)=(x,y)$.

Inversion provides inverses: we have to construct thin homotopies
\begin{equation*}
\gamma \pcomp \prev\gamma \sim \id_y 
\quand 
\id_x \sim \prev\gamma \pcomp \gamma
\end{equation*}
for any path $\gamma$ with $\ev(\gamma)=(x,y)$. To do so, consider a \emph{smoothing function} $\varphi$, i.e. a smooth map $\varphi\maps  [0,1] \to [0,1]$ such that there exists $\varepsilon>0$ with $\varphi(t)=0$ for $t<\varepsilon$ and $\varphi(t)=1$ for $t>1-\varepsilon$. Notice that for $s \in [0,1]$,
$\Gamma_s(t) \df  \gamma(\varphi(s)\varphi(t))$
is a path in $X$ with $\ev(\Gamma_s) = (x,\gamma_{\varphi}(s))$, where $\gamma_{\varphi} \df  \gamma \circ \varphi$ is the reparameterized path. Furthermore, it can be regarded as  a path
\begin{equation*}
\Gamma\maps  [0,1] \to PX\maps  s \mapsto \Gamma_s
\end{equation*}
with $\ev(\Gamma)=(\id_{x},\gamma_{\varphi})$. The adjoint map $\exd\Gamma$ has rank one since it factors through $[0,1]$.  
Then, $c(r(\Gamma),\Gamma)$ in the notation introduced above is a thin homotopy between $\id_{x}=\id_x \pcomp \id_x$ and $\prev{\gamma_{\varphi}}\pcomp \gamma_{\varphi}$. Since $\gamma$ and $\gamma_{\varphi}$  are thin homotopy equivalent, the latter is thin homotopy equivalent to $\prev{\gamma}\pcomp \gamma$. Thus, we have constructed one of the claimed thin homotopies. The other one can be constructed analogously.
\end{proof}

Finally we remark that the space $\pt X$ of thin homotopy classes of paths in $X$ is functorial in $X$: for $f\maps X \to Y$  a smooth map, the induced map $Pf\maps PX \to PY$ on path spaces descends to thin homotopy classes to a smooth map
\begin{equation*}
\pt f\maps  \pt X \to \pt Y\text{,}
\end{equation*}
respecting composition, reversal and identity paths. An alternative formulation of Proposition \ref{prop:compdesc} is that $\pt X$ is the space of morphisms of a diffeological groupoid called the \emph{path groupoid of} $X$ (see \cite{schreiber3,schreiber5}). From this point of view, the maps $\pt f$ define functors between these groupoids .

\subsection{Fusion Maps}

\label{sec:fusion}

In this section we give the definition of a fusion map,  based on  a  relationship between path spaces and loop spaces that we introduce first. 
A \emph{loop} in  $X$ is a smooth map $\tau\maps  S^1 \to X$. Loops in $X$ form the diffeological space 
\begin{equation*}
LX\df D^\infty(S^1,X)\text{.}
\end{equation*}
The following definition generalizes Barrett's notion of thin homotopies \cite{barret1} from smooth manifolds to diffeological spaces.

\begin{definition}
A \emph{homotopy} between loops $\tau_1$ and $\tau_2$ in $X$ is a path $h\in PLX$ with $\ev(h) = (\tau_1,\tau_2)$. A homotopy $h$ is called \emph{thin}, if the adjoint  map
\begin{equation*}
\exd h\maps  [0,1] \times S^1 \to X\maps  (t,z) \mapsto h(t)(z)
\end{equation*}
has rank one.
\end{definition}

For any smooth map $f\maps S^1 \to S^1$ that is homotopic to the identity -- in particular any rotation and any orientation-preserving reparameterization -- there exists a thin homotopy between $\tau \circ f$ and $\tau$. Analogously to Lemma \ref{lem:thinhomeq} one can show that being thin homotopic is an equivalence relation on the diffeological space $LX$ of loops in $X$. 
We denote this equivalence relation by $\sim$, and the diffeological space of thin homotopy classes of loops by
\begin{equation*}
\lt X \df  LX / \sim
\end{equation*}
and called the \emph{thin loop space} of $X$. 

Now we come to the afore-mentioned relationship between  path spaces and loop spaces.
Denoting by $\pcl X\subset PX$ the subspace of closed paths,  a smooth map $cl\maps  \pcl X \to LX$
is obtained by performing a gluing construction similar to the one from the proof of Lemma \ref{lem:thinhomeq}.
Consider a pair $(\gamma_1,\gamma_2)$ of paths with a common initial and a common end point; such pairs of paths form the fibre product $PX^{[2]} \df  PX \times_{X\times X} PX$, taken along the evaluation map $\ev\maps PX \to X \times X$. We have a map $se \maps  PX^{[2]} \to \pcl X$ that takes the pair $(\gamma_1,\gamma_2)$ to the closed path $\prev{\gamma_2} \pcomp \gamma_1$, and is smooth since path composition and inversion are smooth \cite[V.3, V.4]{iglesias1}. All together, we have a smooth map
\begin{equation*}
l \df  cl\circ se \maps  PX^{[2]} \to LX\text{.}
\end{equation*}
Since every loop in the image of $l$ is  constant in a neighborhood of $1 \in S^1$, it is clear that $l$ is not surjective. It induces, however, a subduction on  thin homotopy classes:

\begin{lemma}
\label{lem:ell}
There is a  subduction $\ell\maps  \pt X^{[2]} \to \lt X$ such that the diagram
\begin{equation*}
\alxydim{@=1.2cm}{PX^{[2]} \ar[r]^{l} \ar[d]_{\mathrm{pr}^2}  & LX \ar[d]^{\mathrm{pr}} \\ \pt X^{[2]} \ar[r]_{\ell} & \lt X}
\end{equation*}
is commutative. In particular, $\ell$ is surjective. Moreover, $\ell$ satisfies the relations
\begin{equation*}
\ell(\gamma_1,\kappa \pcomp \gamma_2) = \ell(\prev\kappa \pcomp \gamma_1,\gamma_2)
\quand
\ell(\gamma_1,\gamma_2 \pcomp \beta) = \ell(\gamma_1 \pcomp \prev\beta,\gamma_2)
\end{equation*}
for all  possible $\gamma_1,\gamma_2,\kappa,\beta\in \pt X$.
\end{lemma}

\begin{proof}
We show that $\mathrm{pr} \circ l$ descends to the claimed map $\ell$; this makes the diagram commutative. With Proposition \ref{prop:compdesc} it remains to prove that $cl$ descends. Indeed: the  composition of a thin homotopy $h\in P\pcl X$ with $cl$ defines a thin homotopy $cl \circ h \in PLX$.
The first relations follow since already $se(\gamma_1,\kappa \pcomp \gamma_2) = se(\prev\kappa \pcomp \gamma_1,\gamma_2)$. The second follows since the loops $cl(\prev\beta \pcomp \prev{\gamma_2} \pcomp \gamma_1)$ and $cl(\prev{\gamma_2} \pcomp \gamma_1 \circ \prev\beta)$ are related by a rotation of $\frac{2\pi}{3}$, and hence thin homotopic. The statement that $\ell$ is a subduction is  not  needed in this article and can hence be left as an exercise. 
\begin{comment}

Let $c\maps U \to \lt X$ be a plot. We have to show that $c$ lifts locally through plots of $\pt X^{[2]}$ along $\ell$. By definition of the diffeology on $\lt X$, every point $x \in U$ has an open neighborhood $U_x\subset U$ and a plot $\tilde c\maps  U_x \to LX$ such that $\mathrm{pr} \circ \tilde c = c|_{U_x}$. Let $\varphi$ be a smoothing map (see the proof of Proposition \ref{prop:compdesc}). Consider the  map
\begin{equation*}
s_{\varphi}\maps LX \to PX^{[2]}\maps  \tau \mapsto (\gamma_1,\gamma_2)
\end{equation*} 
with
\begin{equation*}
\textstyle
\gamma_1(t) \df  \tau(\frac{1}{2}\varphi(t))
\quand
\gamma_2(t) \df  \tau(1-\frac{1}{2}\varphi(t))
\end{equation*}
under the identification $S^1 \cong \R/\Z$.
This is evidently a smooth map, and so is its composition with the projection $\mathrm{pr}^2\maps  PX^{[2]} \to \pt X^{[2]}$. All together, we have constructed a plot 
\begin{equation*}
c' \df  \mathrm{pr}^2 \circ s_{\varphi} \circ \tilde c\maps  U_x \to \pt X^{[2]}
\end{equation*}
of $\pt X^{[2]}$. It remains to check that $\ell \circ c' = c|_{U_x}$. This amounts to construct a thin homotopy $h$ between a loop $\tau\in LX$ and the loop $\ell(s_{\varphi}(\tau))$. Such a thin homotopy can be constructed from a reparameterization obtained from a smooth homotopy between $\varphi$ and the identity $\id_{[0,1]}$.
\end{comment}
\end{proof}

The subduction $\ell$ is needed to define fusion maps. Let $G$ be a Lie group, and suppose $f\maps  \lt X \to G$ is a smooth map. We introduce the notation
\begin{equation*}
f_{\ell} \df   f \circ \ell\maps  \pt X^{[2]} \to G
\end{equation*}
in order to simplify the following formulae.

\begin{definition}
A smooth map $f\maps  \lt X \to G$ is called \emph{fusion}, if 
\begin{equation*}
f_{\ell}(\gamma_1,\gamma_2) \cdot f_{\ell}(\gamma_2,\gamma_3) = f_{\ell}(\gamma_1,\gamma_3)
\end{equation*}
for all $(\gamma_1,\gamma_2,\gamma_3) \in \pt X^{[3]}$, i.e. for all triples of thin homotopy classes of paths with a common initial point and a common end point.
\end{definition}

It is straightforward to deduce the following properties of fusion maps.
\begin{lemma}
\label{lem:fusionprop}
Let $f\maps  \lt X \to G$ be a fusion map. Then,
\begin{enumerate}
\item[(a)]
$f_{\ell}(\gamma_1,\gamma_2) = f_{\ell}(\gamma_2,\gamma_1)^{-1}$ for all $(\gamma_1,\gamma_2) \in \pt X^{[2]}$.

\item[(b)] 
$f_{\ell}(\gamma,\gamma)=1$ for all $\gamma\in \pt X$.

\end{enumerate}
\end{lemma}

Fusion maps form a subspace  of the diffeological space $D^\infty(\lt X,G)$ of all smooth maps, which we write as $\fus G {\lt X}$. It appears on the right hand side of the isomorphism of Theorem \ref{th:main}.
\begin{definition}
A \emph{fusion homotopy} between fusion maps $f_0$ and $f_1$ is a path $h$ in $\fus G{\lt X}$ with $\ev(h)=(f_0,f_1)$.
\end{definition}

Due to the composition and reversal of paths, fusion homotopies define an equivalence relation $\sim$ on the space of fusion maps. We denote the space of equivalence classes by
\begin{equation*}
\fushom G {\lt X} \df  \fus G {\lt X} / \sim\text{.}
\end{equation*}
It appears on the right hand side of the isomorphism of Theorem \ref{th:man}.

There are two particular situations. The first is when $G$ is replaced by an abelian Lie group $A$. Then, fusion maps $\fus A {\lt X}$  form a group by point-wise multiplication, and a subgroup of the group $D^\infty (\lt X,A)$. The group structure is preserved under fusion homotopies; hence,   $\fushom A {\lt X}$ is also a group.

The second situation is that of a smooth manifold $M$. Then one can express the condition that a map $f\maps  \lt M \to G$ is smooth in terms of the Fréchet manifold structure on $LM$. Indeed, $f$ is smooth if and only if  
$f \circ \mathrm{pr}\maps  LM \to G$
is smooth in the Fréchet sense, where $\mathrm{pr}\maps  LM \to \lt M$ is the projection to thin homotopy classes (see Lemmata \ref{lem:quotmap} and \ref{lem:mapspacediff}).  
In the same way, the smoothness of a fusion homotopy $h \in P \fus G{\lt M}$ can be characterized by saying that the pullback of the adjoint map  $\exd h\maps  [0,1] \times \lt M \to G$ to $[0,1] \times LM$ is smooth in the Fréchet sense. Thus, the sets $\fus G {\lt M}$ and $\fushom G {\lt M}$ have a description in terms of   Fréchet manifolds.

\setsecnumdepth{2}

\section{Diffeological Principal Bundles with Connection}

\label{sec:diffbun}

In this section we introduce diffeological principal bundles with connection, which appear on the left hand side of the bijection of Theorem \ref{th:main}.  In addition, we prove some results that we need in order to prove Theorems \ref{th:main} and \ref{th:man}.

\subsection{Diffeological Principal Bundles}

\label{sec:diffbundles}

Here we define diffeological principal bundles and  show that they form a sheaf of groupoids over diffeological spaces (Theorem \ref{th:diffbunsheaf}). This sheaf is monoidal if the structure group is abelian (Theorem \ref{th:diffbunmonsheaf}). We   show that a diffeological bundle over a \emph{smooth manifold} is the same as a \emph{smooth} principal bundle (Theorem \ref{th:equivalencesmooth}). 

Let $G$ be a Lie group  and let $X$  be a diffeological space.

\begin{definition}
\label{def:diffbun}
A \emph{diffeological principal $G$-bundle} over $X$ is a subduction $p\maps  P \to X$ together with a fibre-preserving right action of $G$ on $P$ such that
\begin{equation*}
\tau\maps  P \times G \to P \times_X P\maps  (p,g) \mapsto (p,pg)
\end{equation*}
is a diffeomorphism.
\end{definition}

The condition on the map $\tau$ ensures that the action is smooth, free and fibrewise transitive. For instance, it determines a smooth map
\begin{equation}
\label{cantriv}
g_P\maps  P^{[2]} \to G\maps  (p_1,p_2) \mapsto \mathrm{pr}_2(\tau^{-1}(p_1,p_2))\text{,}
\end{equation}
whose result is the unique group element $g\in G$ with $p_2=p_1g$. The condition that the projection $p$ be a subduction ensures that $X$ is diffeomorphic to the quotient of $P$ by the group action \cite[I.50]{iglesias1}. 

\begin{remark}
Definition \ref{def:diffbun} coincides with the restriction of \cite[Definition 3.3.1]{iglesias2} from diffeological groups to ordinary Lie groups. In order to see this, it is important to notice that the projection of a diffeological principal $G$-bundle is necessarily a \quot{strong subduction}, the diffeological analogue of a submersion.
\end{remark}

The morphisms between diffeological principal $G$-bundles over $X$ are  $G$-equivariant smooth maps that respect the projections to $X$. For $P_1$ and $P_2$ diffeological principal $G$-bundles over $X$, we claim that every morphism $\varphi\maps P_1 \to P_2$ is invertible. To see this, consider the smooth map
\begin{equation*}
P_2 \times_X P_1 \to P_1\maps  (p_2,p_1) \mapsto p_1g_{P_2}(\varphi(p_1),p_2)\text{.}
\end{equation*}
It satisfies the gluing condition for the subduction $\mathrm{pr}_1\maps P_2 \times_X P_1 \to P_2$, and hence descends to a smooth map $\varphi^{-1}\maps P_2 \to P_1$. This is an inverse of $\varphi$. Thus, diffeological principal $G$-bundles over a diffeological space   $X$ form a groupoid that we denote by $\diffbun G X$. 

If $G$ is replaced by an abelian Lie group $A$, the sets of morphisms has the following structure.

\begin{lemma}
\label{lem:homtorsor}
Let $P_1$ and $P_2$ be diffeological principal $A$-bundles over $X$. Then the set $\hom(P_1,P_2)$ of morphisms is a torsor over the group $D^\infty(X,A)$. 
\end{lemma}

\begin{proof}
The proof goes exactly as in the manifold case and uses the smoothness of the map $\tau$ from Definition \ref{def:diffbun} and the smoothness of the map $g_P$ from \erf{cantriv}. 
\end{proof}

Pullbacks of diffeological principal $G$-bundles are defined in the same way as for smooth principal $G$-bundles. Thus, diffeological principal $G$-bundles form a presheaf of groupoids over diffeological spaces. We shall see that this presheaf is actually a sheaf with respect to the Grothendieck topology of subductions (see Appendix \ref{sec:diffcov}).

The gluing axiom is formulated as follows.
Associated to any subduction $\omega\maps  W \to X$ is a \emph{descent category} $\mathrm{Des}(\omega)$. The objects of $\mathrm{Des}(\omega)$ are pairs $(P,d)$ of a diffeological principal $G$-bundle over $W$ and of a morphism
$d\maps  \omega_1^{*}P \to \omega_2^{*}P$
of diffeological principal $G$-bundles over $W^{[2]}$ such that 
\begin{equation}
\label{desccond}
\omega_{23}^{*}d \circ \omega_{12}^{*}d = \omega_{13}^{*}d
\end{equation}
over $W^{[3]}$. Here, $\omega_{i_1,...,i_k}$ is the projection to the indexed factors. A morphism in $\mathrm{Des}(\omega)$ between objects $(P_1,d_1)$ and $(P_2,d_2)$ is a morphism $\varphi\maps  P_1 \to P_2$ of diffeological principal $G$-bundles over $W$ such that 
\begin{equation}
\label{eq:desmorph}
d_2 \circ\omega_1^{*}\varphi  = \omega_2^{*}\varphi\circ d_1\text{.}
\end{equation}
The pullback along $\omega$ defines a functor 
\begin{equation*}
\omega^{*}\maps  \diffbun G X \to \mathrm{Des}(\omega)\text{.}
\end{equation*}
The gluing axiom is

\begin{lemma}
\label{sheafproof}
For every subduction $\omega\maps W \to X$, the functor $\omega^{*}$
is an equivalence of groupoids. 
\end{lemma}

\begin{proof}
We  construct an inverse functor $\omega_{*}$.
For $(P,d)$ an object in $\mathrm{Des}(\omega)$, we consider
\begin{equation}
\label{eq:eqreldesc}
P' \df  P / \sim
\quad\text{ with }\quad
p_1 \sim p_2 \Leftrightarrow d(p_1) = (p_2)\text{.}
\end{equation}
Due to \erf{desccond}, $\sim$ is an equivalence relation, and $P'$ is equipped with the pushforward diffeology. The projection $\omega \circ p\maps  P \to X$ satisfies the gluing condition for the subduction $\mathrm{pr}\maps  P \to P'$, so that it descends to a smooth map $p'\maps P' \to X$.
The action of $G$ descends to $P'$ since $d$ is $G$-equivariant. Now we have to show that $P'$ is a diffeological principal $G$-bundle over $X$. 

First we show that $p'$ is again a subduction. Let $c\maps  U \to X$ be a plot and $x \in U$. Since $\omega \circ p$ is -- as a composition of subductions -- a subduction, there exists an open neighborhood $V\subset U$ of $x$ and a plot $\tilde c\maps  V \to P$ of $P$. Now, $c'\df  \mathrm{pr} \circ \tilde c$ is a plot of $P'$, and $p' \circ c' = c|_V$. Thus, $p'$ is a subduction. 

In order to verify that the map $\tau'$ associated to $P'$ is a diffeomorphism, consider the commutative diagram
\begin{equation*}
\alxydim{@=1.2cm}{P \times G \ar[r]^-{\tau} \ar[d]_{\mathrm{pr} \times \id} & P \times_W P \ar[d]^{\mathrm{pr} \times \mathrm{pr}} \\ P' \times G \ar[r]_-{\tau'} & P' \times_X P'\text{.}}
\end{equation*}
The vertical maps are subductions. Thus, $\tau'$ is smooth by Lemma \ref{lem:quotmap}. Since $\tau$ is a bijection,  $\tau'$ has to be a bijection, and again by Lemma \ref{lem:quotmap}, the inverse of $\tau'$ is a smooth map. 

Summarizing, $\omega_{*}(P,d) \df  P'$ is a diffeological principal $G$-bundle over $X$. Now let $\varphi\maps  (P_1,d_1) \to (P_2,d_2)$ be a morphism in $\mathrm{Des}(\omega)$. Due to \erf{eq:desmorph}, there exists a unique map $\varphi'\maps P_1' \to P_2'$ such that the diagram
\begin{equation*}
\alxydim{@=1.2cm}{P_1 \ar[r]^{\varphi} \ar[d]_{p_1'} & P_2 \ar[d]^{p_2'}\\ P_1' \ar[r]_{\varphi'} & P_2'}
\end{equation*} 
is commutative. Again, Lemma \ref{lem:quotmap} shows that $\varphi'$ is smooth. It is also $G$-equivariant, and thus a morphism $\omega_{*}(\varphi) \df  \varphi'$ of diffeological principal $G$-bundles over $X$.

What remains is to define natural equivalences $\omega_{*} \circ \omega^{*} \cong \id$ and $\omega^{*} \circ \omega_{*} \cong \id$. Suppose first $(P,d)$ is an object in $\mathrm{Des}(\omega)$. With $\omega^{*}P' = W \times_X P'$, the map 
\begin{equation*}
\xi_{(P,d)}\maps  P \to \omega^{*}P'\maps  x \mapsto (p(x),\mathrm{pr}(x))
\end{equation*}
is smooth and $G$-equivariant, and  natural in $(P,d)$. 
Suppose secondly that $P$ is a diffeological principal $G$-bundle over $X$. The equivalence relation $\sim$ from \erf{eq:eqreldesc} on $\omega^{*}P=W \times_X P$ identifies $(w,x)$ and $(w',x')$ if and only if $\omega(w)=\omega(w')$ and $x=x'$. Consider the projection $\mathrm{pr}_2\maps  \omega^{*}P \to P$, which is $G$-equivariant. It respects the equivalence relation $\sim$ and defines hence a smooth map $\zeta_P\maps  (\omega^{*}P)' \to P$. This map is natural in $P$.
\end{proof}

Summarizing the above results, we have:

\begin{theorem}
\label{th:diffbunsheaf}
Let $G$ be a Lie group. 
The assignment $X \mapsto \diffbun GX$ defines a sheaf of groupoids over the site of diffeological spaces.
\end{theorem}

For an abelian Lie group $A$ there is more structure: the groupoid $\diffbun A X$ of diffeological principal $A$-bundles is monoidal. The usual definition of tensor products of abelian principal bundles carries over to the diffeological context. If $P_1$ and $P_2$ are diffeological principal $A$-bundles over $X$, then
\begin{equation*}
P_1 \otimes P_2 \df  (P_1 \times_X P_2)/\sim
\quad\text{ with }\quad
(p_1.a,p_2) \sim (p_1,p_2.a)\text{,}
\end{equation*}
 equipped with its canonical diffeology according to Example \ref{ex:constructions} (d) and (e), is again a diffeological principal $A$-bundle over $X$. Verifying that  the functor $\omega_{*}$ and the natural equivalences $\omega_{*} \circ \omega^{*} \cong \id$ and $\omega^{*} \circ \omega_{*} \cong \id$ constructed in the proof of Lemma \ref{sheafproof} are monoidal, we have

\begin{theorem}
\label{th:diffbunmonsheaf}
Let $A$ be an abelian Lie group. The assignment $X \mapsto \diffbun AX$ defines a sheaf of monoidal groupoids over the site of diffeological spaces.
\end{theorem}

In the remainder of this section we consider diffeological principal $G$-bundles over a smooth manifold $M$. The functor $\man \to \diff$ from the category of smooth manifolds to the category of diffeological spaces (see Section \ref{sec:diffdefs}) induces a\  functor
\begin{equation*}
\mathcal{D}_M\maps \bun G M \to \diffbun G M\text{.}
\end{equation*}
One only has to notice that the projection of a smooth principal bundle is a subduction (Lemma \ref{lem:smoothsubduction}).

\begin{theorem}
\label{th:equivalencesmooth}
The functor $\mathcal{D}_M$
is an isomorphism of groupoids. 
\end{theorem}

\begin{proof}
 Suppose $p\maps P \to M$ is a diffeological principal $G$-bundle over $M$. We  equip the total space $P$ with a smooth manifold structure. By Lemma \ref{lem:smoothsubduction} we can choose an open cover $\left \lbrace U_{\alpha}\right \rbrace$ of $M$ together with diffeological sections $s_{\alpha}\maps U_{\alpha} \to P$. They induce bijections 
\begin{equation*}
\sigma_{\alpha}\maps  U_{\alpha} \times G \to p^{-1}(U_{\alpha})\maps  (x,g) \mapsto s_{\alpha}(x)g
\end{equation*}
and thus equip each subset $p^{-1}(U_{\alpha}) \subset P$ with a smooth manifold structure. The transition function is
\begin{equation*}
\sigma_{\beta}^{-1} \circ \sigma_{\alpha}\;\maps \; (U_{\alpha} \cap U_{\beta}) \times G \;\to\; (U_{\alpha} \cap U_{\beta}) \times G
\;\maps \;
(x,g)\;\mapsto\; (x,g_P(s_{\beta}(x),s_{\alpha}(x))g)
\end{equation*}
and thus smooth. Hence, the smooth manifold structures on the sets $p^{-1}(U_{\alpha})$ glue together. The same argument shows that they are independent of the choice of the open sets $U_{\alpha}$ and the sections $s_{\alpha}$. We claim that the original diffeology on $P$ coincides with the smooth diffeology induced by the smooth manifold structure we have just defined. Given that claim, the projection $p\maps P \to X$, the local sections $s_{\alpha}\maps U_{\alpha} \to P$, and the action of $G$ on $P$ are smooth. Thus, $p\maps P \to M$ is a smooth principal $G$-bundle over $M$. Similarly, every morphism $\varphi\maps P_1 \to P_2$ of diffeological principal $G$-bundles is smooth. This yields a functor
\begin{equation*}
\mathcal{D}_M^{-1}\maps \diffbun G M \to \bun G M\text{.}
\end{equation*}

We have to show that the two functors $\mathcal{D}_M$ and $\mathcal{D}_M^{-1}$ are strict inverses of each other. One part is exactly the above claim. In order to prove the claim, suppose $P$ is a diffeological principal $G$-bundle. We have to show that a map $c\maps U \to P$ is a plot of $P$ if and only if it is smooth with respect to the smooth manifold structure on $P$ defined above. Suppose first $c$ is a plot. Define the open sets $V_{\alpha} \df  c^{-1}p^{-1}(U_{\alpha})$ that cover $U$, and consider the composite 
\begin{equation}
\label{compo}
\sigma_{\alpha}^{-1} \circ c|_{V_{\alpha}}\maps  V_{\alpha} \to U_{\alpha} \times G\text{.}
\end{equation}  
Since $\sigma_{\alpha}^{-1}$ is a smooth map, \erf{compo} is smooth. Because $\sigma_{\alpha}$ is a chart of the smooth manifold  $P$, $c$ is smooth on $V_{\alpha}$. Since $U$ is covered by the sets $V_{\alpha}$, $c$ is smooth everywhere. Conversely, suppose $c\maps U \to P$ is smooth. Then, \erf{compo} is smooth and thus a plot of $U_{\alpha} \times G$. But since $\sigma_{\alpha}$  is also smooth, its composition with the plot \erf{compo} is a plot of $P$. Since this composition is $c|_{V_{\alpha}}$, $c$ is a plot by axiom (D3) of Definition \ref{def:diff}.

It remains to check the other part. 
Assume  that $P$ is a smooth principal $G$-bundle. Then, the sections $s_{\alpha}$ used above can be chosen smooth, resulting in diffeomorphisms $\sigma_{\alpha}$. These induce the original smooth manifold structure on $P$. 
\end{proof}

Since pullbacks (and in the abelian case: tensor products) of diffeological principal bundles are defined exactly as for smooth bundles, it is clear that the isomorphisms $\mathcal{D}_M$ define an isomorphism of sheaves of (monoidal) groupoids.

\subsection{Connections, Parallel Transport and Holonomy}

\label{sec:connections}

In this section we introduce connections on diffeological principal bundles, and generalize the statements of Section \ref{sec:diffbundles} to a setup with connections (Theorems \ref{th:diffbunconsheaf}, \ref{th:diffbunconsheafmon} and \ref{th:conequivalencesmooth}). Further we investigate parallel transport and holonomy in diffeological principal bundles with connection.

The definition of a connection is literally the same as in the context of smooth manifolds.

\begin{definition}
\label{def:con}
Let $p\maps  P \to X$ be a diffeological principal $G$-bundle over $X$. A \emph{connection} on $P$ is a 1-form $\omega \in \Omega^1(P,\mathfrak{g})$ such that
\begin{equation*}
\rho^{*}\omega = \mathrm{Ad}^{-1}_{g}(\mathrm{pr}^{*}\omega) + g^{*}\theta\text{,}
\end{equation*}
where $\rho\maps  P \times G \to P$ is the action,  $g\maps P \times G \to G$ and $\mathrm{pr}\maps P\times G \to P$ are the projections, $\mathrm{Ad}$ denotes the adjoint action of $G$ on $\mathfrak{g}$, and $\theta \in \Omega^1(G,\mathfrak{g})$ is the left-invariant Maurer-Cartan form on $G$.
\end{definition}

 Just like this definition, several statements generalize straightforwardly from connections on smooth  principal bundles to diffeological ones. For instance, the trivial principal $G$-bundle $P \df  X \times G$ over $X$ carries a canonical connection $\omega \df  \mathrm{pr}_2^{*}\theta$. More generally, if $A \in \Omega^1(X,\mathfrak{g})$ is any 1-form, 
\begin{equation*}
\omega \df  \mathrm{Ad}^{-1}_{\mathrm{pr}_2}(\mathrm{pr_1}^{*}A) + \mathrm{pr}_2^{*}\theta
\end{equation*}
defines a connection on $P$. 

Let $P_1$ and $P_2$ be principal $G$-bundles with connections $\omega_1$ and $\omega_2$, respectively.
A bundle morphism $\varphi\maps P_1 \to P_2$ is \emph{connection-preserving} if $\varphi^{*}\omega_2 = \omega_1$.
We denote the groupoid of diffeological principal $G$-bundles with connection by  $\diffbuncon G X$. We have the following extension of Theorem \ref{th:diffbunsheaf}.

\begin{theorem}
\label{th:diffbunconsheaf}
Let $G$ be a Lie group. 
The assignment $X \mapsto \diffbuncon GX$ defines a sheaf of groupoids over the site of diffeological spaces.
\end{theorem}

\begin{proof}
Pullbacks of connections are defined in the evident way.
It remains to verify the gluing axiom. Let $\pi\maps Y \to X$ be a subduction, let $P$ be a principal $G$-bundle over $Y$ and let $\omega$ be a connection on $P$. Suppose $d\maps \pi_1^{*}P \to \pi_2^{*}P$ is a connection-preserving bundle morphism satisfying the cocycle condition
\begin{equation*}
\pi_{23}^{*}d \circ \pi_{12}^{*}d = \pi_{13}^{*}d\text{.}
\end{equation*}
Let $P'$ the quotient principal $G$-bundle over $X$, coming with a subduction $\mathrm{pr}\maps P \to P'$. For $\mathrm{pr}_i\maps  P \times_{P'}P \to P$ the two projections, we have to show that $\mathrm{pr}_1^{*}\omega = \mathrm{pr}_2^{*}\omega$. Then, since differential forms form a sheaf \cite[VI.38]{iglesias1}, the 1-form $\omega$ descends to  $P'$. It follows then automatically that the quotient 1-form is a connection.

In order to prove the identity $\mathrm{pr}_1^{*}\omega = \mathrm{pr}_2^{*}\omega$, consider the smooth map 
\begin{equation*}
k\maps P \times_{P'} P \to \mathrm{pr}_1^{*}P\maps (x_1,x_2) \mapsto (x_1,p(x_1),p(x_2))\text{,}
\end{equation*}
where the projections $\mathrm{pr}_{1},\mathrm{pr}_{2}\maps  P \times_{P'} P \to P$ are given by $\mathrm{pr}_1 \circ k$ and $\mathrm{pr}_1 \circ d \circ k$, respectively. Since $d$ preserves connections, we have $\mathrm{pr}_1^{*}\omega = \mathrm{pr}_2^{*}\omega$.
\end{proof}

Next come some statements about connections on diffeological principal bundles with abelian structure group $A$. First of all, it is straightforward to verify that on a tensor product $P_1 \otimes P_2$ of two such bundles, one has a tensor product connection $\omega_1 \otimes \omega_2$, coming from the sum $\mathrm{pr}_1^{*}\omega_1 + \mathrm{pr}_2^{*}\omega_2$ that descends along the subduction $P_1 \times_X P_2 \to P_1 \otimes P_2$. We have immediately

\begin{theorem}
\label{th:diffbunconsheafmon}
The groupoid $\diffbuncon AX$ is monoidal, and principal $A$-bundles with connection form a sheaf of monoidal groupoids.
\end{theorem}

Further, we have the following generalization of Lemma \ref{lem:homtorsor}.

\begin{lemma}
Let $P_1$ and $P_2$ be principal $A$-bundles over $X$  with connections. Then, the set $\hom(P_1,P_2)$ of connection-preserving morphisms is a torsor over the group $D^\infty\lc(X,A)$ of locally constant smooth maps. 
\end{lemma}

\begin{proof}
Given Lemma \ref{lem:homtorsor}, it is enough to prove the following claim. Suppose $\varphi\maps P_1 \to P_2$ is a connection-preserving bundle morphism, and  $f\maps X \to A$ is a smooth map. One computes that $\varphi f$ is connection-preserving if and only if
\begin{equation}
\label{eq:pullbackconn}
p_1^{*}f^{*}\theta = 0\text{,}
\end{equation}
where $p_1\maps P_1 \to X$ is the bundle projection. Because differential forms form a sheaf over  $\diff$, and $p_1$ is a subduction, it follows that \erf{eq:pullbackconn} holds if and only if already $f^{*}\theta = 0$.  According to the following lemma, this is the case if and only if $f$ is locally constant.
\end{proof}

Generally, we define $\mathrm{dlog}(f) \df  f^{*}\theta$ for any smooth function $f\maps  X \to A$. Then:

\begin{lemma}
\label{lem:locconst}
 $\mathrm{dlog}(f) = 0$ if and only if $f$ is locally constant.
\end{lemma}

\begin{proof}
Suppose first that $f$ is locally constant. Then, for any plot $c\maps U \to X$, its pullback $f \circ c$ is constant on path-connected components of $U$, i.e. locally constant. Thus, the ordinary differential form $(f^{*}\theta)_c = (f \circ c)^{*}\theta$ vanishes. Conversely, suppose $f^{*}\theta=0$. Assume that there exists a path $\gamma \in PX$ with $\ev(\gamma)=(x,y)$ such that $f(x) \neq f(y)$. It follows that the composition $\tau \df  f \circ \gamma$ is a smooth, non-constant map. In particular, there exists $t \in (0,1)$ and $v\in T_t(0,1)$ such that $\mathrm{d}\tau|_t(v)\neq 0$. Then,
\begin{equation*}
(f^{*}\theta)_{\gamma}|_t(v) = \tau^{*}\theta|_t(v) \neq 0\text{,}
\end{equation*}
contradicting the assumption of $f^{*}\theta=0$.
\end{proof}

Next we return to a general Lie group $G$, and compare diffeological principal bundles with connection over a smooth manifold to smooth bundles with connection. The functor $\mathcal{D}_M$ from Section \ref{sec:diffbundles} extends to a functor
\begin{equation*}
\mathcal{D}_M^{\nabla}\maps  \buncon G M \to \diffbuncon G M\text{,}
\end{equation*}
and as a consequence of Theorem \ref{th:equivalencesmooth} we see immediately

\begin{theorem}
\label{th:conequivalencesmooth}
The functor $\mathcal{D}_M^{\nabla}$ is an isomorphism of groupoids.
\end{theorem}

The \emph{curvature} of a connection $\omega$ on a principal $G$-bundle $P$ is the 2-form
\begin{equation*}
K_\omega \df  \mathrm{d}\omega + [\omega \wedge \omega] \in \Omega^2(P,\mathfrak{g})\text{.}
\end{equation*}
If $G$ is abelian $K_{\omega}$ descends to a 2-form $\Omega^2(X,\mathfrak{g})$.  A connection $\omega$ is called \emph{flat} if $K_{\omega}$
vanishes. Flatness can be detected \quot{locally}:

\begin{lemma}
\label{lem:flatloc}
A connection $\omega$ on a diffeological principal $G$-bundle $p\maps P\to X$ is flat if and only if for every plot $c\maps U \to X$ the pullback connection $c^{*}\omega$ on  $c^{*}P$ is flat.
\end{lemma}

\begin{proof}
Clearly, if $\omega$ is flat, $c^{*}\omega$ is flat. Let $d\maps V \to P$ be a plot of $P$, so that $c \df  p \circ d$ is a plot of $X$. By assumption, $c^{*}P$ is flat; additionally  it also has a smooth section $s\maps  V \to c^{*}P\maps  v \mapsto (v,d(v))$. Then,
\begin{equation*}
0= s^{*}K_{c^{*}\omega} = s^{*}\mathrm{pr}^{*}K_{\omega} = d^{*}K_{\omega} = (K_{\omega})_d\text{.}
\end{equation*}
This shows that the 2-form $K_{\omega}$ vanishes.
\end{proof}

Important examples of flat connections arise as follows.

\begin{lemma}
\label{lem:homflat}
Let $f\maps X \to Y$ be a smooth rank one map, and $P$ a principal $G$-bundle over $Y$ with connection. Then, $f^{*}P$ is flat.
\end{lemma}

\begin{proof}
Using Lemma \ref{lem:flatloc} we may check that $c^{*}f^{*}P$ is flat for all plots $c\maps  U \to X$. Moreover, since $c^{*}f^{*}P$ is a smooth principal $G$-bundle with connection (Theorem \ref{th:conequivalencesmooth}), we can check its flatness locally. Since $f$ has rank one, every point $u \in U$ has an open neighborhood $V \subset U$ such that $(f \circ c)|_V$ factors through a rank one map $g\maps  V \to W$ and a plot $d\maps W \to Y$ of $Y$. It follows that $c^{*}f^{*}P|_V \cong g^{*}d^{*}P$, which is flat. 
\end{proof}

In the remainder of this section we define parallel transport and holonomy for connections on diffeological principal bundles. For this purpose we regard a connection $\omega \in \Omega^1(P,\mathfrak{g})$ on a diffeological principal $G$-bundle $P$ over $X$ via Theorem \ref{th:formsfunctors} as a smooth map $F_{\omega}\maps  \pt P \to G$. The condition on the 1-form $\omega$ from Definition \ref{def:con} is now saying that for  $g \in PG$ and  $\gamma \in PP$ we have
\begin{equation}
\label{eq:confun}
g(1) \cdot F_{\omega}(\gamma g) = F_{\omega}(\gamma)\cdot g(0)\text{.}
\end{equation}

In order to define the parallel transport, let us first notice the following general fact. If $f\maps M \to X$ is a smooth map defined on a contractible smooth manifold $M$, then $f$ lifts to $P$, i.e. there exists a smooth map $\tilde f\maps  M \to P$ such that $p \circ \tilde f=f$. Indeed, the pullback $f^{*}P$ is by Theorem \ref{th:equivalencesmooth} a smooth principal $G$-bundle and thus has a smooth section $s\maps  M \to f^{*}P$. Combining this section with the projection $\mathrm{pr}\maps  f^{*}P \to P$ yields the claimed lift.

\begin{definition}
\label{def:ptr}
Suppose $P$ is a diffeological principal $G$-bundle with connection $\omega$. Let $\gamma\in PX$ be a path and $\tilde\gamma$ be a lift. Then, the map
\begin{equation*}
\ptrcon\gamma\omega\maps  P_{\gamma(0)} \to P_{\gamma(1)}\maps  q \mapsto \tilde\gamma(1).(F_{\omega}(\tilde\gamma)\cdot g_P(\tilde\gamma(0),q))
\end{equation*}
is called the \emph{parallel transport} of $\omega$ along $\gamma$. 
\end{definition}

It is straightforward to check that the parallel transport $\ptrcon\gamma\omega$ is independent of the choice of the lift $\tilde\gamma$. Indeed, if $\gamma'$ is another lift, we have a smooth map
\begin{equation*}
g\maps  [0,1] \to G\maps  t \mapsto g_P(\tilde\gamma(t),\gamma'(t))
\end{equation*}
such that $\tilde\gamma g = \gamma'$ and thus
\begin{eqnarray*}
\gamma'(1).(F_{\omega}(\gamma')\cdot g_P(\gamma'(0),q)) &=& \tilde\gamma g(1).(F_{\omega}(\tilde\gamma g)\cdot g_P(\tilde\gamma g(0),q))
\\
&=& \tilde\gamma(1).(g(1) \cdot F_{\omega}(\tilde\gamma g)\cdot g(0) \cdot g_P(\tilde\gamma (0),q))
\\
&=& \tilde\gamma(1).(F_{\omega}(\tilde\gamma) \cdot g_P(\tilde\gamma (0),q))
\end{eqnarray*} 
with the last equality given by \erf{eq:confun}.

Alternatively, parallel transport can be defined using the  parallel transport of smooth principal bundles. For this purpose, one pulls back $(P,\omega)$ along $\gamma$ to a principal $G$-bundle over $[0,1]$. Denote by
\begin{equation*}
\ptrcon\tau{\gamma^{*}\omega}\maps  (\gamma^{*}P)|_0 \to (\gamma^{*}P)|_1
\end{equation*}
the parallel transport of $\gamma^{*}\omega$ along the canonical path $\tau$ from $0$ to $1$. Then, under the identification $(\gamma^{*}P)|_t \cong P_{\gamma(t)}$, we have $\ptrcon\gamma\omega = \ptrcon\tau{\gamma^{*}\omega}$.

We summarize all properties of parallel transport in the following  

\begin{proposition}
\label{prop:proppar}
Let $P$ be a diffeological principal $G$-bundle over $X$ with connection $\omega$. Then,
\begin{enumerate}
\item[(a)]
parallel transport is functorial in the path: 
\begin{equation*}
\ptrcon{\id_x}\omega = \id_{P_x}
\quand
\ptrcon{\gamma_2}\omega \circ \ptrcon{\gamma_1}\omega = \ptrcon{\gamma_2 \circ \gamma_1}\omega
\end{equation*}

\item[(c)]
the map $\ptrcon\gamma\omega$ is a $G$-equivariant diffeomorphism, and  depends only on the thin homotopy class of the path $\gamma$.

\item[(c)]
if the connection $\omega$ is flat, $\ptrcon{\gamma}\omega$ only depends on the homotopy class of $\gamma$.

\item[(d)]
the map
\begin{equation*}
\ptrcon{}\omega\maps  \pt X \times_X P \to P \maps  (\gamma,q) \mapsto \ptrcon\gamma\omega(q)
\end{equation*}
is smooth. 

\end{enumerate}
\end{proposition}

\begin{proof}
(a) follows directly from the functorial properties of the map $F_{\omega}$. Equivariance and invertibility in (b) are clear from the definition. 

To see the independence from thin homotopies, consider a thin homotopy $h$ between paths $\gamma_1$ and $\gamma_2$ and its adjoint $\exd h\maps  [0,1]^2 \to X$ that we extend to a plot $c\maps  \R^2 \to X$ (see the proof of Theorem \ref{th:formsfunctors}). As noticed above, one can choose a smooth section $s\maps  \R^2 \to c^{*}P$ and obtain a lift $\tilde c \df  \mathrm{pr} \circ s\maps  \R^2 \to P$ of $c$. The lift is a homotopy between lifts $\tilde\gamma_1$ and $\tilde\gamma_2$ (though not thin, in general). We claim that $\tilde c^{*}F_{\omega}$ is flat. Then, by Lemma \ref{lem:flathom}, $F_{\omega}(\tilde\gamma_1)=F_{\omega}(\tilde\gamma_2)$. To prove the claim, we notice that the connection $\mathrm{pr}^{*}\omega$ on  $c^{*}P$ is flat by Lemma \ref{lem:homflat}. But then,  $\tilde c^{*}F_{\omega} = s^{*}\mathrm{pr}^{*}\omega$ is also flat.
For (c) the same proof applies, just that $\tilde c^{*}F_{\omega}$ is already flat by assumption.

Finally, to see (d) consider the  map
\begin{equation*}
\tilde\tau\maps  PP \times_X P \to P\maps  (\tilde\gamma,q)\mapsto \tilde\gamma(1).(F_{\omega}(\tilde\gamma)\cdot g_P(\tilde\gamma(0),q))
\end{equation*} 
which is  smooth as a composition of smooth maps. Furthermore, for $\gamma\in PX$ and any lift $\tilde\gamma\in PP$, we have by definition $\ptrcon\gamma\omega(q) = \tilde\tau(\tilde\gamma,q)$.  Then we claim that $P\pi\maps  PP \to PX$ is a subduction. Since $\tilde\tau$ is independent of the lift $\tilde\gamma$, it descends by Lemma \ref{lem:quotmap} to a smooth map. To see that $P\pi$ is a subduction, notice that any plot $c\maps U \to X$ lifts locally over contractible open neighborhoods to $P$, as noticed previously.
\end{proof}

The following discussion is restricted to diffeological principal $A$-bundles with connection for $A$ an abelian Lie group.
In order to define the holonomy of a connection we have to regard loops as closed paths. Since a path has by definition sitting instants, one has to choose a smoothing function $\varphi$; then we obtain a smooth map
\begin{equation*}
o_{\varphi}\maps  LX \to \pcl X \maps  \tau \mapsto \tau \circ \varphi\text{.}
\end{equation*}
The thin homotopy class of $o_{\varphi}(\tau)$ is independent of the choice of $\varphi$. Thus we have a canonical smooth map
\begin{equation*}
o\maps  LX \to \pt X\text{.}
\end{equation*} 
Notice that $o$ does not descend to the thin loop space $\lt X$, since a rotation of the loop would change the endpoints of the associated path.

\begin{definition}
\label{def:hol}
Let $\tau\maps S^1 \to X$ be a loop, and let $q\in P_{\tau(0)}$ be an element in the fibre of $P$ over the base point of $\tau$. The holonomy of the connection $\omega$ around $\tau$ is the unique group element $\mathrm{Hol}_{\omega}(\tau) \df a \in A$ such that
\begin{equation*}
\ptrcon{o(\tau)}\omega(q).a=q\text{.}
\end{equation*}
\end{definition}

Since the structure group $A$ is abelian, the holonomy is independent of the choice of $q$ (we recall that for non-abelian groups the holonomy is only  well-defined up to a conjugation). Analogously to parallel transport, it can be expressed in terms of the holonomy of the smooth principal $A$-bundle $\tau^{*}P$ over $S^1$, namely
$\mathrm{Hol}_{\omega}(\tau) = \mathrm{Hol}_{\tau^{*}\omega}(S^1)$.
The following proposition summarizes further important properties of the holonomy.

\begin{proposition}
\label{prop:holsmooth}
Let $P$ be a diffeological principal $A$-bundle and $\omega$ a connection on $P$. The holonomy $\mathrm{Hol}_{\omega}(\tau)$ around a loop $\tau$  depends only  on the thin homotopy class of $\tau$, and defines a smooth map
\begin{equation*}
\mathrm{Hol}_{\omega}\maps  \lt X \to A\text{.}
\end{equation*}
Furthermore,
\begin{enumerate}

\item[(a)]
If $(P_1,\omega_1)$ and $(P_2,\omega_2)$ are isomorphic as principal $A$-bundles with connection,
\begin{equation*}
\mathrm{Hol}_{\omega_1} = \mathrm{Hol}_{\omega_2}\text{.}
\end{equation*}

\item[(b)]
If $(P_1,\omega_1)$ and $(P_2,\omega_2)$ are principal $A$-bundles with connections,
\begin{equation*}
\mathrm{Hol}_{\omega_1\otimes \omega_2} = \mathrm{Hol}_{\omega_1} \cdot \mathrm{Hol}_{\omega_2}\text{.}
\end{equation*}

\end{enumerate}
\end{proposition}

\begin{proof}
We first verify that $\mathrm{Hol}_{\omega}\maps LX \to A$ is smooth. We have to show that for every plot $c\maps  U \to LX$ the map $\mathrm{Hol}_{\omega} \circ c\maps  U \to A$ is smooth. This can be checked locally. For $u\in U$, let $V \subset U$ be a contractible open neighborhood of $u$. Then, the base point projection $\ev_0 \circ c\maps V \to X$ lifts to a smooth map $\tilde q\maps  V \to P$. Now,
\begin{equation*}
(\mathrm{Hol}_{\omega} \circ c)|_V(v) = \ptrcon{}\omega(o(c(v)),\tilde q(v))
\end{equation*}
is a composition of smooth maps (see Proposition \ref{prop:proppar} (d)), and hence smooth. Now we verify that $\mathrm{Hol}_{\omega}(\tau)$ depends only on the thin homotopy class of $\tau$. For $h\in PLX$ a thin homotopy between loops $\tau_1$ and $\tau_2$, and $\exd h\maps  [0,1] \times S^1 \to X$ its adjoint, we see by Lemma \ref{lem:homflat} that $\exd h^{*}P$ is flat. Hence, by Stokes' Theorem, 
\begin{equation*}
1 = \mathrm{Hol}_{\exd h^{*}\omega}(S^1 \times \left \lbrace 0 \right \rbrace)^{-1} \cdot \mathrm{Hol}_{\exd h^{*}\omega}(S^1 \times \left \lbrace 1 \right \rbrace)\text{.}
\end{equation*}
Thus, $\mathrm{Hol}_{\omega}(\tau_1) = \mathrm{Hol}_{\omega}(\tau_2)$. The algebraic properties (a) and (b) are straightforward to check.
\end{proof}

Our final observation concerns the \quot{derivative} of holonomy in the general diffeological setup, generalizing a well-known formula in the manifold setup.
\begin{proposition}
\label{prop:holder}
Let $P$ be a principal $A$-bundle with connection $\omega$ over a diffeological space $X$. Then,
\begin{equation*}
\mathrm{dlog}(\mathrm{Hol}_{\omega}) = \int_{S^1} \ev^{*}K_{\omega} \;\in\; \Omega^1(\lt X,\mathfrak{a})\text{.}
\end{equation*}
In particular, $\mathrm{Hol}_{\omega}$ is locally constant if $\omega$ is flat. 
\end{proposition}

\begin{proof}
We work in a plot $d\maps  U \to \lt X$. Let $\gamma\maps  (-\varepsilon,\varepsilon) \to U$ represent a tangent vector in $U$, with $u \df  \gamma(0)$.
Let $C_{0,\varepsilon} \df  [0,t] \times S^1$ be the standard cylinder, with $\tau_t$ denoting the loop at time $t$. On $C_{0,\varepsilon}$ we pick the orientation which makes $\tau_{\varepsilon}$ orientation-preserving.
Consider the map $\phi\maps  C_{0,\varepsilon} \to X$ defined by $\phi(t,z) \df  \tilde d(\gamma(t),z)$, where $\tilde d\maps  U \times S^1 \to X$ is the smooth map associated to $d$. Then, Stokes' Theorem (applied to the \emph{smooth} principal $A$-bundle $\phi^{*}P$ over $C_{0,\varepsilon}$) yields
\begin{equation*}
\mathrm{Hol}_{\omega}(\tau_t) = \mathrm{Hol}_{\omega}(\tau_0) \cdot \exp \left (  \int_{C_{0,\varepsilon}} \phi^{*}K_{\omega} \right )\text{.}
\end{equation*}
Now we compute
\begin{equation*}
\left .\mathrm{dlog}(\mathrm{Hol}_{\omega})_d \right |_u \left ( \left .\frac{\mathrm{d}}{\mathrm{d}t} \right|_0 \gamma \right) = \theta_{\mathrm{Hol}_{\omega}(\tau_0)} \left ( \left . \frac{\mathrm{d}}{\mathrm{d}t} \right|_0 \mathrm{Hol}_{\omega}(\tau_t) \right )= \left . \frac{\mathrm{d}}{\mathrm{d}t} \right|_0  \int_{C_{0,t}} \phi^{*}K_{\omega}  \text{,}
\end{equation*}
where the first step is the definition of $\mathrm{dlog}$, and the second step is Stokes' Theorem. The result is precisely the fibre integration of $\tilde d^{*}K_{\omega}$, evaluated at $u$ and on the tangent vector given by $\gamma$. Thus, for each plot $d$ we have $\mathrm{dlog}(\mathrm{Hol}_{\omega})_d = \int_{S^1} \tilde d^{*}K_{\omega}$, and this shows the claim. The statement about flatness follows from Lemma \ref{lem:locconst}. 
\end{proof}

\setsecnumdepth{2}

\section{Regression}

\label{sec:regression}

Throughout this section, $X$ is a connected diffeological space, $x$ is a base point, and $G$ is a Lie group. In the first subsection, we construct a diffeological principal $G$-bundle $\un_x(f)$ over $X$ associated to a fusion map $f\maps  \lt X \to G$.  In the second subsection, we equip this bundle with a connection. These two constructions constitute the bijections of Theorems \ref{th:man} and \ref{th:main}, respectively.

\subsection{Reconstruction of the Bundle}

\label{sec:recon}

The reconstruction of the bundle $\un_x(f)$ is essentially the one of  \cite{barret1}, recalled in a way emphasizing the role of descent theory. 
Let  $\ptx Xx$ denote the subspace of $\pt X$ consisting of classes of paths in $X$ starting at $x$. Then, the restriction of the subduction $\ev\maps \pt X \to X \times X$ to $\ptx Xx$ is still a subduction 
$\ev_1\maps  \ptx Xx \to X$,
which is a consequence of Proposition \ref{prop:grothendieck}. 

Let $f\maps  \lt X \to G$ be a fusion map. Let $T_x \df  \ptx Xx \times G$ denote the trivial principal $G$-bundle over $\ptx Xx$. We equip $T_x$ with a descent structure for the subduction $\ev_1$, and use that diffeological principal $G$-bundles form a sheaf over diffeological spaces (Theorem \ref{th:diffbunsheaf}).  The descent structure is a bundle morphism
\begin{equation}
\label{eq:descstruc}
d_{f}\maps  \mathrm{pr}_1^{*}T_x \to \mathrm{pr}_2^{*}T_x
\end{equation}   
over $\ptx Xx^{[2]}$, with $\mathrm{pr}_i\maps  \ptx Xx^{[2]} \to \ptx Xx$ the projections. It is defined by
\begin{equation*}
d_f(\gamma_1,\gamma_2,g) \df  (\gamma_1,\gamma_2,f_{\ell}(\gamma_2,\gamma_1)g)\text{.}
\end{equation*}
This is a smooth, $G$-equivariant map, respects  the projections to the base, and satisfies the descent condition \erf{desccond} due to the fusion property of $f$. Thus, we obtain a diffeological  principal $G$-bundle 
\begin{equation*}
\un_x(f) \df  (\ev_1)_{*}(T_x,d_{f})
\end{equation*}
over $X$.

\begin{lemma}
\label{lem:basepointind}
The isomorphism class of $\un_x(f)$ does not depend on the choice of the base point $x$.
\end{lemma}

\begin{proof}
For another base point $y$ choose a path $\kappa \in \pt X$ with $\ev(\kappa) = (y,x)$, which is possible since $X$, is by assumption, connected. Proposition \ref{prop:compdesc} shows that the map
\begin{equation*}
c_{\kappa}\maps  \ptx Xx \to \ptx Xy\maps  \gamma \mapsto \gamma \pcomp \kappa
\end{equation*}
is smooth. It suffices to show that the descent structure $d_f$ on the trivial principal $G$-bundle $T_y$ over $\ptx Xy$ pulls back  along $c_{\kappa}$ to the descent structure $d_f$ on $T_x$. Indeed, this is equivalent to the identity
$f_{\ell}(\gamma_2 \pcomp \kappa,\gamma_1 \pcomp \kappa) = f_{\ell}(\gamma_2,\gamma_1)$
which follows from Lemma \ref{lem:ell}.
\end{proof}

Next we look at the two particular cases we have looked at at the end of Section \ref{sec:fusion}. The first case is that $G$ is replaced by an abelian Lie group $A$. We recall that then the product of fusion maps is again a fusion map.

\begin{lemma}
\label{lem:reconmon}
Let $f_1,f_2\maps  \lt X \to A$ be fusion maps. There is a canonical isomorphism
\begin{equation*}
\un_x(f_1f_2) \cong \un_x(f_1) \otimes \un_x(f_2)\text{.}
\end{equation*}
\end{lemma}

\begin{proof}
The tensor product $T_x \otimes T_x$ inherits a descent structure $d_{f_1} \otimes d_{f_2}$, and since diffeological bundles form a sheaf of \emph{monoidal} groupoids (Theorem \ref{th:diffbunmonsheaf}), $T_x \otimes T_x$ descends to $\un_x(f_1) \otimes \un_x(f_2)$. Consider the isomorphism 
\begin{equation*}
\varphi\maps  T_x \otimes T_x \to T_x\maps  ((\gamma,a_1),(\gamma,a_2)) \mapsto (\gamma,a_1a_2)
\end{equation*}
of  principal $A$-bundles over $\ptx Xx$.  It exchanges the descent structure $d_{f_1} \otimes d_{f_2}$ on $T_x \otimes T_x$ with $d_{f_1f_2}$ on $T_x$; hence, $\varphi$  descends to the claimed isomorphism.
\end{proof}

The second particular case is when $X$ is a smooth manifold $M$. Then, there is an isomorphism $\diffbun G M \cong \bun G M$ between the groupoids of diffeological principal $G$-bundles and ordinary smooth ones (Theorem \ref{th:equivalencesmooth}). The existence of connections on smooth principal bundles permits us to show the following:
\begin{lemma}
Let $M$ be a connected smooth manifold, and let $f\maps  \lt M \to G$ be a fusion map. Then, the isomorphism class of $\un_x(f)$ depends only on the fusion homotopy class of $f$.
\end{lemma}

\begin{proof}
Let $h$ be a fusion homotopy between fusion maps $f_0$ and $f_1$. Let $T$ denote the trivial principal $G$-bundle over $[0,1] \times \ptx Mx$. We define a descent structure on $T$ with respect to the subduction $\id \times \ev_1\maps  [0,1] \times \ptx Mx \to [0,1] \times M$ by
\begin{equation*}
d_{h}(t,\gamma_1,\gamma_2,g) \df  (t,\gamma_1,\gamma_2,h(t)(\ell(\gamma_2,\gamma_1)) g)\text{.}
\end{equation*}
The cocycle condition is satisfied because $h(t)$ is a fusion map at any time $t$. We obtain a diffeological (and thus smooth) principal $G$-bundle $Q \df  (\id \times \ev_1)_{*}(T,d_h)$ over $M$. Its restriction to $\left \lbrace 0 \right \rbrace \times M$ is $\un_x(f_0)$ and its restriction to $\left \lbrace 1\right \rbrace \times M$ is $\un_x(f_1)$. Thus, any connection on $Q$ defines an isomorphism.
\end{proof}

Summarizing, on a smooth manifold $M$ we have a well-defined map
\begin{equation*}
\un\maps  \fushom G {\lt M} \to \hc 0\bun G M\text{.}
\end{equation*}
If $G$ is abelian,  $\un$ is a group homomorphism by  Lemma \ref{lem:reconmon}. This group homomorphism defines the bijection of Theorem \ref{th:man}.

\subsection{Reconstruction of the Connection}

\label{sec:reconcon}

Our construction of a connection on $\un_x (f)$ is different from the one of \cite{barret1}.
We equip the trivial $G$-bundle $T_x$ over $\ptx Xx$ with a connection. This amounts to specifying a 1-form $A_{x} \in \Omega^1(\ptx Xx,\mathfrak{g})$. Then we prove that this connection descends to a connection on $\un_x(f)$.

\begin{figure}
\begin{equation*}
\includegraphics{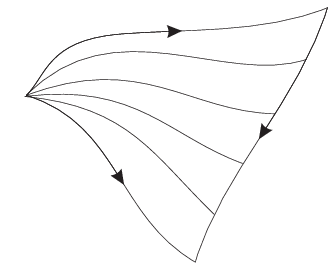}\setlength{\unitlength}{1pt}\begin{picture}(0,0)(211,681)\put(50.43770,758.61326){$x$}\put(116.46170,796.10995){$\beta_0$}\put(91.46863,714.45061){$\beta_1$}\put(173.90778,731.17729){$\beta$}\end{picture}
\end{equation*}
\caption{A path in the space $\ptx Xx$ of based paths.}
\label{fig:pfadpfad}
\end{figure}

In order to define the 1-form we use Theorem \ref{th:formsfunctors}, which identifies 1-forms on any diffeological space with certain smooth maps on its path space. Let $\gamma\maps [0,1] \to \ptx Xx$ be a path. Notice that $\beta_0 \df  \gamma(0)$ and $\beta_1 \df  \gamma(1)$ are elements in $\ptx Xx$. Further notice that
$\beta(t) \df  \gamma(t)(1)$
defines a  path $\beta \in PX$ (see Figure \ref{fig:pfadpfad}). 
We consider the group element
\begin{equation*}
F_x(\gamma) \df  f_{\ell}(\beta_1,\beta \pcomp \beta_0) \in G\text{.}
\end{equation*}

\begin{lemma}
\label{lem:Fxdef}
This defines a smooth map $F_x\maps  \pt (\ptx Xx) \to G$ satisfying 
\begin{equation*}
F_{x}(\gamma' \circ \gamma) = F_{x}(\gamma') \cdot F_{x}(\gamma)
\end{equation*}
for every pair $(\gamma,\gamma')$ of composable elements in $\pt(\ptx Xx)$. 
\end{lemma}

\begin{proof}
The map
\begin{equation}
\label{eq:impmap}
P(\ptx Xx) \to \ptx Xx \times\ptx Xx \times PX\maps  \gamma \mapsto (\beta_0,\beta_1,\beta) \end{equation}
we have implicitly used above
is smooth. Smoothness of the composition \cite[V.3]{iglesias1} and of $f_{\ell}$ show then that $F_x$ is a smooth map. 
From Lemmata \ref{lem:ell} and \ref{lem:fusionprop} we deduce
\begin{equation*}
F_x(\gamma' \pcomp \gamma) = f_{\ell}(\beta_2,\beta' \pcomp \beta \pcomp \beta_0) = f_{\ell}(\beta^{\prime-1} \pcomp \beta_2,\beta_1)\cdot f_{\ell}(\beta_1,\beta \pcomp \beta_0)=F_x(\gamma')\cdot F_x(\gamma)\text{.}
\end{equation*}
It remains to show that  $F_x$ is well-defined on $\pt(\ptx Xx)$. 
If $h$ is a thin homotopy between $\gamma$ and $\gamma'$, and  $(\beta_0,\beta_1,\beta)$ and $(\beta_1,\beta_2,\beta')$ denote their images under \erf{eq:impmap}, then $\beta_0=\beta_0'$ and $\beta_1=\beta_1'$. Furthermore, the paths $\beta$ and $\beta'$ are thin homotopic: $h_1(s)(t) \df  \ev_1(h(s)(t))$ defines a homotopy between $\beta$ and $\beta'$, for which $\exd {h_1} = \ev_1 \circ \exd h$ has rank one since $\exd h$ has rank one (Lemma \ref{lem:thinprop} (d)). Thus, $F_x(\gamma)=F_x(\gamma')$.
\end{proof}

By Theorem \ref{th:formsfunctors},  $F_x$ defines a 1-form $A_x \in \Omega^1(\ptx Xx,\mathfrak{g})$. We want to show that the connection $\omega_x$ on $T_x$ determined by $A_{x}$ descends along $\ev_1$. For this purpose, we have to show that
\begin{lemma}
The bundle morphism $d_f\maps  \mathrm{pr}_1^{*}T_x \to \mathrm{pr}_2^{*}T_x$ preserves connections.
\end{lemma}

\begin{proof}
Whenever one has a morphism $P_{1} \to P_2$ between trivial principal $G$-bundles over a diffeological space $Y$, given by multiplication with a smooth map $f\maps Y \to G$, this morphism preserves connection 1-forms $A_1$ and $A_2$ on $P_1$ and $P_2$, respectively, if and only if
\begin{equation*}
A_2 = \mathrm{Ad}^{-1}_{f}(A_1) + f^{*}\theta\text{.}
\end{equation*}
This can be checked explicitly in the same way as one does it in the smooth manifold context. In our situation, we have $Y \df \ptx Xx^{[2]}$, $P_i \df  \mathrm{pr}_i^{*}T_x$,  $f \df  f_{\ell}^{-1}$, and the 1-forms are $A_i \df  \mathrm{pr}_i^{*}A_x$ for $i=1,2$. Thus, the equation we have to show is 
\begin{equation*}
\mathrm{pr}_1^{*}A_{x} = \mathrm{Ad}_{f_{\ell}}(\mathrm{pr}_2^{*}A_{x}) - f_{\ell}^{*}\bar\theta\text{.}
\end{equation*}
By Theorem \ref{th:formsfunctors}, this is equivalent to showing that $F_x$ satisfies
\begin{equation*}
f_{\ell}(\gamma(1),\gamma'(1)) \cdot F_{x}(\gamma') = F_{x}(\gamma) \cdot f_{\ell}(\gamma(0),\gamma'(0))
\end{equation*}
for all $(\gamma,\gamma') \in \pt (\ptx Xx^{[2]})$. Indeed, since then $\beta=\beta'$, we have
\begin{equation*}
f_{\ell}(\beta_1,\beta_1') \cdot f_{\ell}(\beta'_1,\beta' \pcomp \beta'_0) = f_{\ell}(\beta_1,\beta\pcomp\beta_0') =f_{\ell}(\beta^{-1} \pcomp \beta_1,\beta_0')  = f_{\ell}(\beta_1,\beta \pcomp \beta_0) \cdot f_{\ell}(\beta_0,\beta_0')\text{,} \end{equation*}
following from Lemmata \ref{lem:ell} and \ref{lem:fusionprop}.
\end{proof}

Since diffeological principal bundles with connection form a sheaf of groupoids (Theorem \ref{th:diffbunconsheaf}) we have now a principal $G$-bundle 
\begin{equation*}
\uncon_x(f) \df  (\ev_1)_{*}(T_x,\omega_x,d_f)
\end{equation*}
over $X$ \emph{with connection}, whose underlying bundle is $\un_x(f)$ from the previous section. 

\begin{lemma}
The isomorphism class of $\uncon_x(f)$ does not depend on the choice of the base point $x$.
\end{lemma}

\begin{proof}
The map $c_\kappa$ from Lemma \ref{lem:basepointind} satisfies $c_{\kappa}^{*}A_x = A_y$, equivalently, $c_{\kappa}^{*}F_x = F_y$, which follows from Lemma \ref{lem:ell}. Thus, $c_{\kappa}$ descends to the claimed isomorphism.
\end{proof}

Summarizing, we have a well-defined map
\begin{equation*}
\uncon\maps  \fus G {\lt X} \to  \hc 0 \diffbuncon G X
\end{equation*}
that we call \emph{regression}. It defines the bijection of Theorem \ref{th:main}. For a smooth manifold $M$, it fits by construction into the commutative diagram
\begin{equation*}
\alxydim{@=1.3cm}{\fus G {\lt M} \ar[r]^-{\uncon} \ar[d] & \hc 0 \buncon G M \ar[d] \\ \fushom G {\lt M} \ar[r]_-{\un} & \hc 0 \bun G M\text{,}}
\end{equation*}
and thus contributes one part of the proof of Theorem \ref{th:comp}.
The following lemma shows that $\uncon$ is a group homomorphism for abelian Lie groups. 

\begin{lemma}
Let $f_1,f_2\maps  \lt X \to A$ be fusion maps. Then, the canonical isomorphism
\begin{equation*}
\un_x(f_1f_2) \cong \un_x(f_1) \otimes \un_x(f_2)
\end{equation*}
from Lemma \ref{lem:reconmon} respects the connections. 
\end{lemma}

\begin{proof}
We recall from the proof of Lemma \ref{lem:reconmon} that the isomorphism has been obtained from an isomorphism $\varphi\maps  T_x \otimes T_x \to T_x$ over $\ptx Xx$. Since the connection $\omega_x$ on $T_x$ is the pullback of a 1-form $A_x$ on $\ptx Xx$, it follows that
\begin{equation*}
\varphi^{*}\omega_x = \mathrm{pr}_1^{*}\omega_x + \mathrm{pr}_2^{*}\omega_x
\end{equation*} 
whatever the definition of $A_x$ was. Thus, $\varphi$ preserves connections and descends to a connection-preserving isomorphism. 
\end{proof}

Finally, we prove one part  of Corollary \ref{co:flat}. Generally, a  map  on a diffeological space $X$ is called \emph{locally constant}, if $f(x)=f(y)$ whenever $(x,y)$ is in the image of the evaluation map $\ev\maps  P X \to X \times X$. 

\begin{proposition}
\label{prop:flatreg}
Let $f\maps  \lt X \to G$ be a locally constant fusion map. Then, the connection on $\uncon_x(f)$ is flat. 
\end{proposition}

\begin{proof}
It suffices to show that the 1-form $A_x$ on $\ptx Xx$ is flat, i.e. 
$\mathrm{d}A_x + [A_x \wedge A_x] = 0$.
This is, by Theorem \ref{th:formsfunctors}, equivalent to showing that the  smooth map $F_x\maps  \pt\ptx Xx \to G$ takes the same value on  homotopic paths. Suppose $h\in PP \ptx Xx$ is a homotopy between paths $\gamma,\gamma' \in P\ptx Xx$. If $(\beta_0,\beta_1,\beta)$ and $(\beta_0',\beta_1',\beta)$ are the triples of paths associated to $\gamma$ and $\gamma'$,
we find $\beta_0=\beta_0'$ and $\beta_1=\beta_1'$, and $h$ induces a homotopy $\tilde h$ between $\beta$ and $\beta'$ (see the proof of Lemma \ref{lem:Fxdef}). Then,
\begin{equation*}
\alxydim{@C=1.8cm}{[0,1] \ar[r]^{\tilde h} & PX \ar[r]^{\mathrm{pr}} & \pt X \ar[r]^{\ell(\beta_1,-\pcomp \beta_0)} & \lt X}
\end{equation*}
is a path in $\lt X$ from $\ell(\beta_1,\beta \pcomp \beta_0)$ to $\ell(\beta_1',\beta'\pcomp \beta_0')$. Since $f$ is locally constant, it follows that $F_x(\gamma)=F_x(\gamma')$. \end{proof}

\setsecnumdepth{1}

\section{Transgression}

\label{sec:transgression}

In this section, $X$ is a diffeological space and $A$ is an abelian Lie group.
Let $P$ be a principal $A$-bundle over $X$ with connection $\omega$. According to Proposition \ref{prop:holsmooth}, the holonomy of $\omega$ is a smooth map
$\mathrm{Hol}_{\omega}\maps  \lt X \to A$.
Furthermore,  $\mathrm{Hol}_{\omega}$ depends only on the isomorphism class of $(P,\omega)$ and satisfies 
\begin{equation*}
\mathrm{Hol}_{\omega_1 \otimes \omega_2} = \mathrm{Hol}_{\omega_1} \cdot \mathrm{Hol}_{\omega_2}
\end{equation*}
for $(P_1,\omega_2)$ and $(P_2,\omega_2)$ two diffeological principal $A$-bundles with connection.

\begin{lemma}
\label{lem:holfusion}
The holonomy of a connection $\omega$ is a fusion map.
\end{lemma}

\begin{proof}
This is a simply calculation: for $q$ an element in the fibre of $P$ over the common initial point of three paths $(\gamma_1,\gamma_2,\gamma_3)\in \pt X^{[3]}$,
\begin{eqnarray*}
q= \ptrcon{\prev{\gamma_3} \pcomp \gamma_1}\omega(q).\mathrm{Hol}_{\omega}(\ell(\gamma_1,\gamma_3)) &=& \ptrcon{\prev{\gamma_3} \pcomp \gamma_2}\omega (\ptrcon{\prev{\gamma_2}\pcomp\gamma_1}\omega(q)).\mathrm{Hol}_{\omega}(\ell(\gamma_1,\gamma_3)) \\&=& \ptrcon{\prev{\gamma_3} \pcomp \gamma_2}\omega (q).\mathrm{Hol}_{\omega}(\ell(\gamma_1,\gamma_2))^{-1}.\mathrm{Hol}_{\omega}(\ell(\gamma_1,\gamma_3))
\\
&=& q.\mathrm{Hol}_{\omega}(\ell(\gamma_2,\gamma_3))^{-1}.\mathrm{Hol}_{\omega}(\ell(\gamma_1,\gamma_2))^{-1}.\mathrm{Hol}_{\omega}(\ell(\gamma_1,\gamma_3))
\end{eqnarray*}
Here we have used the definition of holonomy (Definition \ref{def:hol}) via the parallel transport $\ptrcon\gamma\omega$ of the connection $\omega$,  its functorality and its $A$-equivariance (Proposition \ref{prop:proppar}).
\end{proof}

Summarizing, we have a well-defined group homomorphism
\begin{equation*}
\trcon\maps  \hc 0 \diffbuncon A X \to \fus A {\lt X}
\end{equation*}
that we call \emph{transgression}. We prove in the following section that it is the inverse of the group homomorphism $\un^{\nabla}$ constructed in Section \ref{sec:reconcon}.

\begin{lemma}
The fusion homotopy class of $\mathrm{Hol}_{\omega}$ is independent of the choice of the connection on $P$.
\end{lemma}

\begin{proof}
For two connections $\omega_0$ and $\omega_1$ on $P$, consider the principal $G$-bundle 
\begin{equation*}
\id \times p\maps  [0,1] \times P \to [0,1] \times X
\end{equation*}
and the 1-form $\Omega \df  t\omega_1 + (1-t)\omega_0$ on $[0,1] \times P$, which defines a connection on $[0,1] \times P$.  Consider further the smooth map
\begin{equation*}
\eta\maps  [0,1] \times \lt X \to \lt([0,1] \times X)
\text{, }\quad
\eta(t,\tau)(s) \df  (\varphi(t),\tau(s))\text{,}
\end{equation*}
where $\varphi$ is a smoothing function (see the proof of Proposition \ref{prop:compdesc}).
Then, 
\begin{equation*}
H \df  \mathrm{Hol}_{\Omega} \circ \eta\maps  [0,1] \times \lt X \to A
\end{equation*}
is a smooth map and corresponds to a path  $h\in P D^\infty(\lt X,A)$ with $H=\exd h$, connecting $\mathrm{Hol}_{\omega_0}$ and $\mathrm{Hol}_{\omega_1}$. All that remains is to check that $h(t)$ lies in the subspace of fusion maps for all $t\in[0,1]$. This follows from the fact that $\mathrm{Hol}_{\Omega}$ is a fusion map (Lemma \ref{lem:holfusion}). 
\end{proof}

Over a smooth manifold $M$, every diffeological principal $A$-bundle is an ordinary, smooth principal bundle (Theorem \ref{th:equivalencesmooth}). Due to the existence of connections on such bundles, we obtain a well-defined group homomorphism
\begin{equation*}
\tr\maps  \hc 0 \bun A M \to \fushom A {\lt M}\text{.}
\end{equation*}
By construction, the diagram
\begin{equation*}
\alxydim{@=1.3cm}{\hc 0 \buncon A M \ar@{->}[r]^-{\trcon} \ar[d] & \fus A {\lt M} \ar[d] \\ \hc 0 \bun A M \ar@{->}[r]_-{\tr} & \fushom A {\lt M}\text{,}}
\end{equation*}
is commutative, which contributes the remaining part to the proof of Theorem \ref{th:comp}.

\setsecnumdepth{1}

\section{Proof of Theorem \ref{th:main}}

\label{sec:proofs}

We show that regression $\uncon$ and transgression $\trcon$ are inverses to each other, starting with the proof that  $\trcon \circ \uncon$ is the identity on the space $\fus A {\lt X}$ of fusion maps. Let $f\maps  \lt X \to A$ be a fusion map. We have to compute the holonomy of the reconstructed bundle $(P,\omega)\df \uncon_x(f)$. Let  a loop $\tau\in L X$
be represented by a closed path $\gamma \in \pcl X$ under the map $cl$ from Section \ref{sec:fusion}. Let $\ptrcon\gamma{\omega}\maps P_{y} \to P_{y}$ denote the parallel transport along $\gamma$, where $y=\gamma(0)=\gamma(1)$ and $P_y$ denotes the fibre of $P$ over $y$. We have to compute $a\in A$ such that $\ptr\gamma(q).a = q$ for some (and hence all) $q \in P_{y}$.

The path $\gamma$ lifts to $\pt X_x$. To see this, let us denote by $\gamma_t \in PX$ the path
$\gamma_t(s) \df  \gamma(t\varphi(s))$,
for $\varphi$ a smoothing function (see the proof of Proposition \ref{prop:compdesc}). We  choose a  path $\kappa\in \pt X$ with $\ev(\kappa)=(x,y)$. Then,
\begin{equation*}
\tilde\gamma\maps  [0,1] \to \pt X_x\maps  t \mapsto \gamma_t \pcomp \kappa
\end{equation*}
is a path and lifts $\gamma$ along the evaluation $\ev_1$. 
We recall that $P$ is descended from the trivial principal $A$-bundle $T_x$ over $\pt X_x$. In particular, it comes with a projection $\mathrm{pr}\maps T_x \to P$.  The parallel transports $\ptrcon\gamma\omega$ in $P$ and $\ptr{\tilde\gamma}$ in $T_x$ fit into a commutative diagram
\begin{equation*}
\alxydim{@=1.3cm}{T_x|_{\kappa} \ar[r]^{\ptr{\tilde\gamma}} \ar[d]_{\mathrm{pr}} & T_x|_{\gamma \circ \kappa} \ar[d]^{\mathrm{pr}} \\ P_{y} \ar[r]_{\ptrcon\gamma\omega} & P_{y}\text{.}}
\end{equation*}
The parallel transport in the trivial bundle $T_x$ along $\tilde\gamma$ is according to Definition \ref{def:ptr} given by 
\begin{equation*}
\ptr{\tilde\gamma}(\tilde\gamma(0),g) \df  (\tilde\gamma(1),gF_x(\tilde\gamma))\text{.}
\end{equation*}
We compute from the definition of $F_x$ and Lemma \ref{lem:fusionprop} that
\begin{equation*}
F_x(\tilde\gamma) = f_{\ell}(\tilde\gamma(1),\gamma \pcomp \tilde\gamma(0)) = f_{\ell}(\gamma \pcomp \kappa,\gamma \pcomp \kappa) = 1\text{.}
\end{equation*}
Thus, for $q\df \mathrm{pr}(\tilde\gamma(0),g)\in P_y$, we have $\ptrcon\gamma\omega(q) = \mathrm{pr}(\tilde\gamma(1),g) \in P_{y}$.
Now we use the definition of $P$, given by the descent structure $d_f$ from \erf{eq:descstruc} and the descent construction \erf{eq:eqreldesc}. We compute using Lemma \ref{lem:ell}
\begin{multline*}
\ptrcon\gamma\omega(q)=\mathrm{pr}(\tilde\gamma(1),g) = \mathrm{pr}(\tilde\gamma(0), f_{\ell}(\tilde\gamma(0),\tilde\gamma(1)) g) \\= q.f_{\ell}(\kappa,\gamma \pcomp \kappa)=q.f_{\ell }(\id,\gamma)=q.f_{\ell}(\gamma,\id)^{-1}=q.f(\tau)^{-1}\text{.}
\end{multline*}
We conclude that $\mathrm{Hol}_{\omega}(\tau)=f(\tau)$,
which completes the proof that $(\trcon \circ \uncon)(f) = f$.

Next is the proof that $\uncon \circ \trcon$ is the identity on the group $\hc 0 \diffbuncon A X$ of isomorphism classes of diffeological principal $A$-bundles over $X$ with connection. Let $P$ be such a bundle, and let $f\maps  \lt X \to A$ be the associated fusion map, its holonomy. Let $q_0\in P$ be a fixed point in the fibre over $x$.
Notice that
\begin{equation*}
\varphi\maps  T_x \to \ev_1^{*}P\maps  (\gamma,g) \mapsto (\gamma,\tau_{\gamma}(q_{0}).g)
\end{equation*}
defines a bundle morphism over $\pt X_x$:  it is smooth by Proposition \ref{prop:proppar} (c) and $A$-equivariant. Moreover, it exchanges the descent structure $d_f$ on the trivial bundle $T_{x} =\pt X_x \times A$ with the trivial descent structure on $\ev_1^{*}P$:
\begin{multline*}
\mathrm{pr}_2^{*}\varphi(d_f(\gamma_1,\gamma_2,g)) = \mathrm{pr}_2^{*}\varphi(\gamma_1,\gamma_2,f_{\ell}(\gamma_2,\gamma_1)g) \\= \tau_{\gamma_2}(q_{0}).f_{\ell}(\gamma_2,\gamma_1)g =  \tau_{\gamma_1}(q_{0}).g = \mathrm{pr}_1^{*}\varphi(\gamma_1,\gamma_2,g)\text{.}
\end{multline*}
Hence, $\varphi$ descends to an isomorphism 
\begin{equation*}
(\ev_1)_{*}(\varphi)\maps  \un_x(f) \to P\text{.}
\end{equation*}

It remains to prove that this isomorphism respects the connections. This is the case if and only if the isomorphism $\varphi$ respects connections, i.e. pulls back the connection  $\ev_1^{*}\omega$ on $\ev_1^{*}P$ to the connection $\omega_x$ on $T_x$, which was defined by a 1-form $A_x$.  Consider the composite 
\begin{equation*}
\alxydim{}{\pt X_x \ar[r] & \pt X_x \times A \ar[r]^-{\varphi} & \mathrm{ev}_1^{*}P \ar[r] & P}
\end{equation*}
in which the first map sends a path $\gamma$ to the pair $(\gamma,1)$. Explicitly, this composite is the map
\begin{equation*}
j\maps  \pt X_x \to P\maps  \gamma \mapsto \tau_{\gamma}(q_{0})\text{.}
\end{equation*}
Now, the isomorphism $\varphi$ is connection-preserving if and only if 
\begin{equation}
\label{eq:connpres}
j^{*}\omega = A_x\text{.}
\end{equation}
In order to check equation \erf{eq:connpres}, let  $F_{\omega}\maps  \pt P \to A$ be the smooth map corresponding to the 1-form $\omega$ under the bijection of Theorem \ref{th:formsfunctors}. We recall that $A_x$ was defined by a smooth map $F_x\maps  \pt \ptx Xx \to A$. We claim that $F_x(\gamma) = F_{\omega}(j \circ \gamma)$
for all  $\gamma\maps  [0,1] \to \pt X_x$. This shows \erf{eq:connpres} and completes the proof of Theorem \ref{th:main}.

To show the claim, we recall from Section \ref{sec:reconcon} that $\gamma$ defines three paths, namely $\beta_0 \df  \gamma(0) \in \pt X_x$, $\beta_1 \df  \gamma(1) \in \pt X_x$ and $\beta \in PX$. The path $j \circ \gamma$ is a lift of the path $\beta$ along the bundle projection $p\maps P \to X$: 
\begin{equation*}
p((j \circ \gamma)(t)) = p(\tau_{\gamma(t)}(q_{0})) = \gamma(t)(1) = \beta(t)
\end{equation*}
for all $t\in[0,1]$.
Using the functorality of parallel transport and the definition of holonomy we find
\begin{equation}
\label{eq:fast}
\ptrcon{\beta}\omega(j(\gamma(0))).\mathrm{Hol}_{\omega}(\beta \pcomp \beta_0 \pcomp \beta_1^{-1}) =\ptrcon{\beta_1 \pcomp \beta_0^{-1}}\omega(j(\gamma(0)))\text{.}
\end{equation}
On the left hand side, 
\begin{equation*}
\mathrm{Hol}_{\omega}(\beta \pcomp \beta_0 \pcomp \beta_1^{-1}) = \mathrm{Hol}_{\omega}(\beta_1^{-1} \pcomp \beta \pcomp \beta_0)  = f_{\ell}(\beta \pcomp \beta_0,\beta_1) = f_{\ell}(\beta_1,\beta \pcomp \beta_0)^{-1}\text{,}
\end{equation*} 
where the first equality holds because the two loops in the argument of $\mathrm{Hol}_{\omega}$ are related by a rotation, i.e. a particular thin homotopy. The second equality is the definition of $f$  and the last equality is Lemma \ref{lem:fusionprop}.
On the right hand side of \erf{eq:fast},  we have
\begin{equation*}
\ptrcon{\beta_1 \circ \beta_0^{-1}}\omega(j(\gamma(0))) =\ptrcon{\beta_1}\omega(q_{0}) = j(\gamma(1)) 
\end{equation*}
due to the definition of $j$. Thus, \erf{eq:fast} becomes $\ptrcon{\beta}\omega(j(\gamma(0))) =j(\gamma(1)).f_{\ell}(\beta_1,\beta \pcomp \beta_0)$. 
A comparison with Definition \ref{def:ptr} shows that
$F_{\omega}(j \circ \gamma) = f_{\ell}(\beta_1,\beta \pcomp \beta_0)$.
The right hand side is precisely the definition of $F_x(\gamma)$ and shows that $F_x(\gamma)=F_{\omega}(j \circ \gamma)$.

\begin{appendix}

\setsecnumdepth{2}

\section{Diffeology of Loop Spaces and Paths Spaces}

\label{sec:diff}

In Section \ref{sec:diffdefs} we review standard definitions and facts about diffeological spaces.
In Section \ref{sec:diffcov} we introduce a Grothendieck topology on the category of diffeological spaces based on subductions. In Section \ref{sec:diffforms} we review differential forms on diffeological spaces and relate them to smooth maps on path spaces. 

\subsection{Diffeological Spaces}

\label{sec:diffdefs}

Diffeological spaces were introduced by Souriau \cite{souriau1}, and are today understood as one flavor of the \quot{convenient calculus}. For a concise presentation we refer the reader to \cite{baez6,laubinger1},  for a comprehensive treatment to \cite{iglesias1}, and to \cite{stacey1} for a comparison with other forms of the convenient calculus.

\begin{definition}
\label{def:diff}
A \emph{diffeological space} is a set $X$ with a diffeology. A \emph{diffeology} on a set $X$ is a set of maps
$c \maps  U \to X$
called \emph{plots}, where each plot is defined on an open subset $U \subset \R^{k}$ with varying $k\in
\N_0$, such that three axioms are satisfied:
\begin{itemize}
\item[(D1)]
for any plot $c\maps  U \to X$, any open subset $V \subset \R^l$, and any smooth function $f\maps V \to U$ the map  
$c \circ f$ is also a plot. 
\item[(D2)]
every constant map $c\maps U \to X$ is a plot.
\item[(D3)]
if $f\maps U \to X$ is a map defined on $U \subset \R^k$ and  $\lbrace U_i \rbrace_{i\in I}$ is an open cover of $U$ for which all restrictions $f|_{U_{i}}$ are plots of $X$, then  $f$ is also a plot.
\end{itemize}
Moreover, a map $f\maps X \to Y$ between diffeological spaces  $X$ and $Y$ is  called  \emph{smooth} if
for every plot $c\maps  U \to X$ of $X$ the map $f \circ c\maps U \to Y$ is a plot of $Y$.  
\end{definition}

Diffeological spaces form a category $\diff$, and the isomorphisms in $\diff$ are called \emph{diffeomorphisms}.
\begin{comment}
A diffeology on a set $X$ determines a topology on $X$, namely the finest topology in which all plots are continuous. Putting these topologies, every diffeological map is continuous. In other words, there is a functor
\begin{equation}
\label{eq:funtop}
\diff \to \top\text{,}
\end{equation} 
with $\top$ the category of topological spaces. 
\end{comment}
\begin{comment}
We remark that the composition of the two functors \erf{eq:funfrech} and \erf{eq:funtop} is the evident functor that forgets the manifold structure and only keeps the underlying topology.
\end{comment}
The following theorem explains \emph{why} it is convenient to use diffeological spaces.

\begin{theorem}[{{\cite[Theorem 9]{baez6}}}]
\label{quasitopos}
The category $\diff$ is a quasitopos. 
\end{theorem}

In practice, this  means that many constructions are available which are in general obstructed or simply not possible for smooth manifolds or Fréchet manifolds. The next example describes some of these constructions. 

\begin{example}
\label{ex:constructions}
\begin{enumerate}[(a)]
\item
For diffeological spaces $X$ and $Y$, the set $D^\infty(X,Y)$ of  smooth maps $f\maps X \to Y$ carries a canonical diffeology called the \emph{functional diffeology} \cite[I.57]{iglesias1}.

A map 
$c\maps  U \to D^\infty(X,Y)$
is a plot if and only if the composite 
\begin{equation*}
\alxydim{@C=1.5cm}{U \times X \ar[r]^-{c \times \id} & D^\infty(X,Y) \times X \ar[r]^-{\mathrm{ev}} & Y}
\end{equation*}
is smooth.

\item
Every subset $Y$ of a diffeological space $X$ carries a canonical diffeology called the \emph{subset diffeology} \cite[I.33]{iglesias1}. 

A map $c\maps  U \to Y$ is a plot of $Y$ if and only if its composition with the inclusion $\iota\maps Y \incl X$ is a plot of $X$.

\item
The direct product $X \times Y$ of diffeological spaces $X$ and $Y$ carries a canonical diffeology called the \emph{product diffeology} \cite[I.55]{iglesias1}.

A map $c\maps  U \to X \times Y$ is a plot of $X \times Y$ if and only if its composition with the projections to $X$ and to $Y$ are plots of $X$ and $Y$, respectively.

\item
For any pair of diffeological maps $f\maps  X \to Z$ and $g\maps Y \to Z$, the fibre product $X \times_Z Y$ is -- as a subset of $X \times Y$ -- a diffeological space. 

\item
For  $X$ a diffeological space, $Y$ a set, and $p\maps X \to Y$ a map, $Y$ carries a canonical diffeology called the \emph{pushforward diffeology} \cite[I.43]{iglesias1}.

A map $c \maps  U \to Y$ is a plot if and only if every point $x \in U$ has an open neighborhood $V \subset U$ such that either $c|_V$ is constant or there exists a plot $\tilde c \maps  V \to X$ of $X$ with $c|_{V} = p \circ \tilde c$.

\end{enumerate}
\end{example}

The pushforward diffeology of Example \ref{ex:constructions} (e) arises frequently in this article, namely when $\sim$ is an equivalence relation on a diffeological space $X$, and $Y\df X/\sim$ is the set of equivalence classes. Then, $Y$ carries the pushforward diffeology induced by the  projection $\mathrm{pr}\maps X \to Y$.
One can easily show

\begin{lemma}[{{\cite[I.51]{iglesias1}}}]
\label{lem:quotmap}
Let $X_1$ and $X_2$ be diffeological spaces, let
 $Y$ be a set and let $p\maps X_1 \to Y$ be a map. Then, a map $f\maps Y \to X_2$ is smooth with respect to the pushforward diffeology on  $Y$ if and only if $f \circ p$ is smooth. 
\end{lemma}

In the remainder of this section we shall relate diffeological spaces to Fréchet manifolds. Let $\frech$ denote the category of Fréchet manifolds. Then, there is a functor
\begin{equation}
\label{eq:funfrech}
\frech \to \diff
\end{equation}
defined as follows.
On objects, it declares on a Fréchet manifold $X$ the \emph{smooth diffeology}. Its plots are all smooth maps $c\maps U \to X$, defined on  open subsets $U \subset \R^k$, for all $k\in \N_0$. On morphisms it is the identity: any smooth map $f\maps X \to Y$  between Fréchet manifolds is diffeological. Indeed, its composition $f \circ c$ with any plot $c\maps  U \to X$ of $X$ is smooth, and thus a plot of $Y$. The following theorem permits the unambiguous usage of the word \quot{smooth}.

\begin{theorem}[{{\cite[Theorem 3.1.1]{losik2}}}]
\label{th:frechet}
The functor \erf{eq:funfrech} is full and faithful, i.e. a map between Fréchet manifolds $X$ and $Y$ is smooth in the manifold sense if and only if it is smooth in the diffeological sense.
\end{theorem}

\begin{remark}
\begin{enumerate}
\item 
Since smooth manifolds form a full subcategory of $\frech$, Theorem \ref{th:frechet} remains true upon substituting \quot{smooth} for \quot{Fréchet}.

\item
If $M$ is a smooth manifold \emph{with boundary}, it still carries the smooth diffeology. However, since our plots are defined on \emph{open} subsets, any map $f\maps M \to X$ that is smooth in the interior of $M$ is already  smooth on $M$. For this article this is negligible: the only manifold with boundary that appears here is the interval $M\df [0,1]$ as the domain of paths. But paths are by definition  constant near the boundary.
\end{enumerate}
\end{remark}

Finally we shall show that the functor \erf{eq:funfrech} identifies the Fréchet manifold $C^\infty(S^1,M)$ with the diffeological space $D^\infty(S^1,M)$ (see Example \ref{ex:constructions} (a)), so that the loop space $LM$ of a smooth manifold has an unambiguous meaning. More generally, we have the following statement.

\begin{lemma}
\label{lem:mapspacediff}
Let $M$ be a smooth manifold and let $K$ be a compact smooth manifold.
The functional diffeology on $D^\infty(K,M)$ and the smooth diffeology on $C^\infty(K,M)$ coincide in the sense that every plot of one is a plot of the other. 
\end{lemma}

\begin{proof}
Let $c\maps  U \to C^{\infty}(K,M)$ be a map. It is a plot of $C^{\infty}(K,M)$ if and only if it is  smooth. It is a plot of $D^\infty(K,M)$ if and only if
\begin{equation}
\label{frechplot}
\alxydim{@C=1.5cm}{U\times K \ar[r]^-{c \times \id} & C^{\infty}(K,M) \times K \ar[r]^-{\mathrm{ev}} & M}
\end{equation}
is smooth. Suppose first that $c$ is smooth. Since the evaluation map is smooth, also \erf{frechplot} is smooth. Hence, every plot of $C^{\infty}(K,M)$ is a plot of $D^\infty(K,M)$. Conversely, assume that \erf{frechplot} is smooth. We want to show that $c$ is smooth. We recall that the Fréchet manifold structure on $C^{\infty}(K,M)$ is the one of the set $\Gamma(K,M \times K)$ of smooth sections in the trivial $M$-bundle over $K$. We also recall that if $\varphi\maps  E \to F$ is a smooth morphism of fibre bundles over $K$, the induced map $\varphi_{*}\maps \Gamma(K,E) \to \Gamma(K,F)$ is smooth \cite[Example 4.4.5]{hamilton1}. Since \erf{frechplot} is smooth, also
\begin{equation*}
\alxydim{@C=1.5cm}{U \times K \ar[r]^-{c \times \id} & C^{\infty}(K,M) \times K \ar[r]^-{\mathrm{ev} \times \id} & M \times K}
\end{equation*}
is a smooth morphism of (trivial) fibre bundles over $K$. It hence induces a smooth map $\tilde c\maps  C^{\infty}(K,U) \to C^{\infty}(K,M)$. Let $i\maps  U \to C^{\infty}(K,U)$ be the inclusion of constant maps, which is  smooth. Hence, the composition $\tilde c \circ i\maps  U \to C^{\infty}(K,M)$ is a smooth map and coincides with $c$. Thus, every plot of $D^\infty(K,M)$ is a plot of $C^{\infty}(K,M)$.
\end{proof}

\subsection{Subductions}

\label{sec:diffcov}

Presheaves can be defined over any category, while 
 the formulation of the gluing axiom, i.e. the definition of a sheaf, requires the choice of a Grothendieck topology. In this section we introduce a Grothendieck topology on the category $\diff$ of diffeological spaces.

\begin{definition}[{{\cite[I.48]{iglesias1}}}]
\label{def:subduction}
A smooth map $\pi\maps Y \to X$ is called \emph{subduction} if the following  condition is satisfied.
For every plot $c\maps  U \to X$ and every $x \in U$ there exists an open neighborhood $V \subset U$ of $x$ and a plot $\tilde c\maps  V \to Y$ such that $\pi \circ \tilde c = c|_V$.
\end{definition}

One can  show that every subduction is surjective, and that an injective subduction is a diffeomorphism \cite[Proposition 1.2.15]{iglesias2}. One can also show that a smooth map $\pi\maps Y \to X$ is a subduction if and only if the diffeology of $X$ is the pushforward diffeology induced by $\pi$ \cite[I.46]{iglesias1}. In particular, all projections $\mathrm{pr}\maps X \to X/\sim$ to spaces of equivalence classes are subductions. 
Further examples of subductions are  projections to a factor in a product or fibre product \cite[I.56]{iglesias1}, and -- as mentioned  in Section \ref{sec:thinhom} -- the endpoint evaluation $\ev\maps  PX \to X \times X$ for $X$ connected \cite[V.6]{iglesias1}. Subductions over smooth manifolds can be characterized in the following way.

\begin{lemma}
\label{lem:smoothsubduction}
Let $M$ be a smooth manifold, let $Y$ be a diffeological space and let $\pi\maps Y \to M$ be a smooth map. Then, $\pi$ is a subduction if and only if every point $p \in M$ has an open neighborhood $W \subset M$ that admits a smooth section $s\maps  W \to Y$.
\end{lemma}

\begin{comment}
\begin{proof}
Assume $\pi$ is a subduction. For given $p \in M$, let $c\maps  U \to W$ be a chart of $M$ with $p \in W \subset M$. By definition of the smooth diffeology, $c$ is a plot of $M$, so that we can apply the condition for subductions to $c$ and the point $x \df  c^{-1}(p)$. Thus, there exists an open neighborhood $V \subset U$ of $x$ and a plot $\tilde c\maps  V \to Y$ such that $\pi \circ \tilde c=c|_V$. For the open subset $W \df  c(V)$ of $M$, the map
$s \df  \tilde c \circ c^{-1}\maps  c(V) \to Y$
is a smooth section. Conversely, assume that $\pi$ has local smooth sections. Let $c\maps  U \to M$ be a plot of $M$ and $x \in U$. Let $W \subset M$ be an open neighborhood of $c(x)$ with a smooth section $s\maps W \to Y$. For the open subset  $V \df  c^{-1}(W)$ of $U$, consider the map $\tilde c \df  s \circ c|_V\maps  V \to Y$. It is a plot of $Y$ since $c|_V$ is a plot of $M$ and $s$ is smooth. \end{proof}
\end{comment}

The next proposition is the main point of this subsection.

\begin{proposition}
\label{prop:grothendieck}
Subductions form a Grothendieck topology on the category $\diff$. 
\end{proposition}

\begin{proof}
The identity $\id_X$ of a diffeological space $X$ is clearly a subduction. The composition of subductions is  a subduction \cite[Proposition 1.2.15]{iglesias2}. Finally, the pullback of a subduction along any diffeological map is a subduction \cite[Proposition 1.4.8]{iglesias2}.
\end{proof}

\subsection{Differential Forms}

\label{sec:diffforms}

Differential forms on diffeological spaces are defined \quot{plot-wise} using the notion of ordinary differential forms on smooth manifolds. 

\begin{definition}[{{\cite[VI. 28]{iglesias1}}}]
\label{def:forms}
Let $X$ be a diffeological space. A  \emph{$k$-form} on $X$ is a family $\left \lbrace \varphi_c \right \rbrace$ of  $k$-forms $\varphi_c \in \Omega^k(U)$ parameterized by plots $c\maps U \to X$, such that $\varphi_{c_1} = f^{*}\varphi_{c_2}$ for every commutative diagram
\begin{equation}
\label{eq:compfamforms}
\alxydim{@C=0.4cm}{U_{c_1} \ar[rr]^-{f} \ar[dr]_{c_1} && U_{c_2} \ar[dl]^{c_2}\\&X&}
\end{equation}
with $c_1$ and $c_2$ plots and $f$ smooth.
\end{definition}

 The set of $k$-forms on a diffeological space $X$ is denoted $\Omega^k(X)$; similarly one defines $k$-forms with values in a vector space $V$, denoted $\Omega^k(X,V)$.
All familiar features of differential forms generalize  from smooth manifolds to diffeological spaces, equipping the family $\Omega^{*}(X)$ with the structure of a differential graded commutative  algebra (dgca).  Further, if $f\maps X \to Y$ is a smooth map between diffeological spaces, there is a pullback
\begin{equation*}
f^{*}\maps  \Omega^k(Y) \to \Omega^k(X)
\quad\text{ with }\quad
(f^{*}\varphi)_{c} \df  \varphi_{f \circ c}
\end{equation*}
for a plot $c$ of $X$, forming a morphism between dgca's. In other words, forms over diffeological spaces form a  presheaf $\sheaf\Omega^{*}$ of dgca's over diffeological spaces. One can show \cite[VI.38]{iglesias1} that this presheaf is even a sheaf.

For a smooth manifold $M$, it is easy to check that the $k$-forms of Definition \ref{def:forms} are the same as ordinary (smooth) $k$-forms on manifolds. 
\begin{comment}
For the compatibility of Definition \ref{def:forms} with smooth manifolds, we note

\begin{lemma}
\label{lem:formsdiffsmooth}
Let $M$ be a smooth manifold, and let $\Omega^k(M)^{\infty}$ denote the space of ordinary differential $k$-forms on the manifold $M$. Then, 
\begin{equation*}
\Omega^k(M)^{\infty} \to \Omega^k(M)\maps  \varphi \mapsto \left \lbrace c^{*}\varphi \right \rbrace
\end{equation*}
defines an isomorphism presheaves of dgca's.
\end{lemma}

Lemma \ref{lem:formsdiffsmooth} allows us  not to distinguish  between  two kinds of differential forms on a smooth manifold.
\end{comment}
The following lemma generalizes a familiar fact from smooth manifolds to diffeological spaces (see Lemma \ref{lem:thinman} and \cite[Lemma 4.2]{schreiber5})

\begin{lemma}
\label{lem:thinflat}
Let $f\maps X \to Y$ be a smooth rank $k$ map. Then the pullback $f^{*}\varphi$ vanishes for all $\varphi\in\Omega^{k+1}(Y)$.
\end{lemma}

\begin{proof}
Let $\varphi\in\Omega^{k+1}(Y)$ and let $c\maps U \to X$ be a plot. We  show that $\varphi_{f \circ c}=0$. Since $f$  has rank $k$, every point $u\in U$ has an open neighborhood $U_u\subset U$ with a plot $d\maps V \to Y$ and a  rank $k$ map $g\maps U_u \to V$ satisfying $d \circ g = f \circ c|_{U_u}$, see Definition \ref{def:thin}. It follows that $\varphi_{f \circ c}|_{U_u} = g^{*}\varphi_{d}=0$. 
\end{proof}

Differential forms can be transgressed to the loop space, and we have used that in Proposition \ref{prop:holder}. Let $\varphi$ be a $k$-form on a diffeological space $X$. Consider a plot $d\maps  U \to LX$ of the loop space, i.e. a map such that the adjoint map $\tilde d\maps  U \times S^1 \to X$ is smooth. Consider the $(k-1)$-form 
\begin{equation*}
\psi_d \df  \int_{S^1} \tilde d^{*}\varphi \in \Omega^{k-1}(U)\text{,}
\end{equation*}
where $\tilde d^{*}\varphi$ is a $k$-form on $U \times S^1$ -- and thus an ordinary differential $k$-form -- and $\int_{S^1}$ denotes \quot{integration along the fibre}. Let $d_1\maps  U_1 \to LX$ and $d_2\maps  U_2 \to LX$ both be plots, and let $f\maps  U_1 \to U_2$ be a smooth map such that $d_2 \circ f = d_1$, then we have $\tilde d_2 \circ  \tilde f =\tilde d_1$, where $\tilde f \df  f \times \id$. It follows that $f^{*}\psi_{d_2}=\psi_{d_1}$. Hence, $\psi \df \left \lbrace \psi_d \right \rbrace$ is a $(k-1)$-form on the diffeological space $LX$. The same procedure works for the thin loop space $\lt X$ instead of $LX$. In both cases, we use the symbolical notation 
\begin{equation*}
\psi \df  \int_{S^1} \ev^{*}\varphi\text{.}
\end{equation*}

\setsecnumdepth{1}

\section{Differential Forms and Smooth Functors}

In this section we discuss a close relation between $1$-forms on a diffeological space $X$ and smooth maps on the path space $\pt X$. For smooth manifolds, this relation has been studied before in \cite{schreiber3,schreiber5}. Let $G$ be a Lie group. On the one hand, consider the 
the following groupoid $\fun XG$. Its objects are  smooth maps $F\maps  \pt X \to G$ satisfying
\begin{equation*}
F(\gamma_2 \pcomp \gamma_1) = F(\gamma_2) \cdot F(\gamma_1)
\end{equation*}
whenever paths $\gamma_1,\gamma_2$ are composable. 
A morphism $g\maps F_1 \to F_2$ is a smooth map $g\maps X \to G$ such that
\begin{equation*}
g(\gamma(1)) \cdot F_1(\gamma) = F_2(\gamma) \cdot g(\gamma(0))\text{.}
\end{equation*}
Composition is multiplication, i.e. $g_2 \circ g_1 \df  g_1g_1$, and the identity morphisms are given by the constant map $g=1$. 
On the other hand, let $\mathcal{Z}_X^1(G)$ be the following groupoid. The objects are 1-forms $A \in \Omega^1(X,\mathfrak{g})$, with $\mathfrak{g}$ the Lie algebra of $G$. A morphism $g\maps A_1 \to A_2$ is a smooth map $g\maps X \to G$ such that
\begin{equation*}
A_{2} = \mathrm{Ad}_g(A_1) - g^{*}\bar\theta\text{,}
\end{equation*}
where $\bar\theta \in \Omega^1(G,\mathfrak{g})$ is the right-invariant Maurer-Cartan form on $G$. Composition and identities are as in $\fun XG$.

Both groupoids are natural in $X$, i.e. if $f\maps X \to Y$ is a smooth map, there are evident pullback functors
\begin{equation}
\label{eq:pullbacks}
f^{*}\maps  \fun YG \to \fun XG
\quand
f^{*}\maps  \mathcal{Z}_Y^1(G) \to \mathcal{Z}_X^1(G)\text{.}
\end{equation}
These pullback functors compose strictly under the composition of smooth maps. In other words, $\fun -G$ and $\mathcal{Z}^1(G)$ are presheaves of groupoids over the category of diffeological spaces. We remark that both presheaves are in general not sheaves.

For a smooth manifold $M$, the groupoids $\fun MG$ and $\mathcal{Z}_M^1(G)$ have been introduced in \cite{schreiber3}. We recall 

\begin{proposition}[{{\cite[Proposition 4.7]{schreiber3}}}]
\label{prop:funformssmooth}
For a smooth manifold $M$ there is an isomorphism of categories 
\begin{equation*}
\alxydim{@C=1.5cm}{\fun MG \ar@<0.3pc>[r]^-{\mathfrak{D}^{\infty}} &  \mathcal{Z}^1_M(G) \ar@<0.3pc>[l]^-{\mathfrak{P}^{\infty}}\text{.}}
\end{equation*}
\end{proposition}

This  isomorphism can be  characterized in terms of  parallel transport.  Suppose $\omega\in\Omega^1(P,\mathfrak{g})$ is a connection on a smooth principal $G$-bundle $P$ over $M$, $\gamma\in PM$ is a path and $\tilde\gamma\in PP$ is a lift of $\gamma$. Then,
\begin{equation*}
\ptrcon\gamma\omega(\tilde\gamma(0)) = \tilde\gamma(1).\mathfrak{P}^{\infty}
(\omega)(\tilde\gamma)\text{,}
\end{equation*}
where $\ptrcon\gamma\omega$ denotes the parallel transport of $\omega$ along $\gamma$, and $\mathfrak{P}^{\infty}(\omega)\maps \mathcal{P}P \to G$ corresponds to the 1-form $\omega$ under Proposition \ref{prop:funformssmooth}.

We  need two  properties of the isomorphism of Proposition \ref{prop:funformssmooth}. Firstly, it commutes with the pullback functors \erf{eq:pullbacks}, i.e. it is an isomorphism of presheaves over $\diff$ \cite[Proposition 1.7]{schreiber5}.
Secondly, let us call a $1$-form $A \in \Omega^1(X,\mathfrak{g})$ \emph{flat}, if the 2-form
\begin{equation*}
K_A \df  \mathrm{d}A + [A \wedge A]
\end{equation*}
vanishes, and let us call an object $F$ in $\fun XG$  \emph{flat}, if $F(\gamma)=F(\gamma')$ whenever $\gamma$ and $\gamma'$ are  homotopic (Definition \ref{def:thinhomotopy}). Then, the flat objects in $\mathcal{Z}^1_M(G)$ correspond precisely to the flat objects in $\fun MG$ \cite[Lemma B.1 (c)]{schreiber3}.

We prove the following generalization of Proposition \ref{prop:funformssmooth} from smooth manifolds to  diffeological spaces. 

\begin{theorem}
\label{th:formsfunctors}
Let $X$ be a diffeological space. Then, there is a  isomorphism of categories 
\begin{equation*}
\alxydim{@C=1.5cm}{\fun XG \ar@<0.3pc>[r]^-{\mathfrak{D}} &  \mathcal{Z}^1_X(G) \ar@<0.3pc>[l]^-{\mathfrak{P}}}
\end{equation*}
that  restricts over every plot $c\maps U \to X$ to the isomorphism of Proposition \ref{prop:funformssmooth}, i.e.
\begin{equation*}
c^{*} \circ \mathfrak{D} = \mathfrak{D}^{\infty} \circ c^{*}
\quand
c^{*} \circ \mathfrak{P} = \mathfrak{P}^{\infty} \circ c^{*}\text{.}
\end{equation*}
\end{theorem}

\begin{proof}
The functor $\mathfrak{D}\maps \fun XG \to \mathcal{Z}^1_X(G)$ is easy to define. If an  object $F$ in $\fun XG$ is given, one has for each plot $c\maps U \to X$ a 1-form $A_c \df  \mathfrak{D}^{\infty}(c^{*}F) \in \Omega^1(U,\mathfrak{g})$. These 1-forms clearly define an object $\left \lbrace A_c\right \rbrace$ in $\mathcal{Z}^1_X(G)$. Moreover, any morphism in $\fun XG$ is automatically a morphism in $\mathcal{Z}_X^1(G)$. This defines the functor $\mathfrak{D}$, and it restricts by construction over each plot to the functor $\mathfrak{D}^{\infty}$. To finish the proof, it thus remains to construct the functor $\mathfrak{P}$ such that is strictly inverse to $\mathfrak{D}$.

\def\ep#1{\widetilde{#1}{}}

Assume $A=\left \lbrace A_c \right \rbrace$ is an object in $\mathcal{Z}_X^1(G)$. A map $F\maps  PX \to G$ is defined as follows. Any path $\gamma\maps  [0,1] \to X$ extends canonically to a plot $\ep\gamma\maps  \R \to X$. Namely, one simply puts $\ep\gamma(t)\df \gamma(0)$ for all $t < 0$ and $\ep\gamma(t)\df \gamma(1)$ for all $t>1$. There is a unique thin homotopy class $\tau \in \pt \R$ of paths in $\R$ with $\ev(\tau)=(0,1)$. Then, we put
\begin{equation*}
F(\gamma) \df  \mathfrak{P}^{\infty}(A_{\ep\gamma})(\tau).
\end{equation*} 
We have to check that this definition yields an object in $\fun XG$. This check consists of the following three parts.

\textit{1.) Compatibility with the path composition}. For the following calculations we use the notation $F_{\ep\gamma} \df  \mathfrak{P}^{\infty}(A_{\ep\gamma})$. If $\gamma_1,\gamma_2\in PX$ are composable paths, consider the two smooth maps 
\begin{equation*}
\iota_1\maps \R \to \R\maps t \mapsto \frac{1}{2}t
\quand
\iota_2\maps \R \to \R\maps t \mapsto \frac{1}{2}+\frac{1}{2}t\text{.}
\end{equation*}
Due to the uniqueness of the path $\tau$,
one has $\tau = \pt\iota_{2}(\tau) \pcomp \pt\iota_1(\tau)$ in $\pt \R$. Then,
\begin{multline*}
F(\gamma_2 \pcomp \gamma_1) = F_{\ep{\gamma_2 \pcomp \gamma_1}}(\tau) = F_{\ep{\gamma_2 \pcomp \gamma_1}}(\pt\iota_{2}(\tau) \pcomp \pt\iota_1(\tau) ) \\=F_{\ep{\gamma_2 \pcomp \gamma_1}}(\pt\iota_{2}(\tau)) \cdot F_{\ep{\gamma_2 \pcomp \gamma_1}}(\pt\iota_1(\tau)) =\iota_{2}^{*}F_{\ep{\gamma_2 \pcomp \gamma_1}}(\tau) \cdot \iota_1^{*}F_{\ep{\gamma_2 \pcomp \gamma_1}}(\tau) \\\stackrel{(*)}{=} F_{\ep{\gamma_2}}(\tau) \cdot F_{\ep{\gamma_1}}(\tau) = F(\gamma_2) \cdot F(\gamma_1)\text{,}
\end{multline*}
where $(*)$ comes from the commutative diagram
\begin{equation*}
\alxydim{@C=0.4cm}{\R \ar[rr]^-{\iota_k} \ar[dr]_{\ep{\gamma_k}} && \R \ar[dl]^{\ep{\gamma_2 \pcomp \gamma_1}}\\&X&}
\end{equation*}
for $k=1,2$, which implies equalities $\iota_k^{*}A_{\ep{\gamma_2\pcomp\gamma_1}}=A_{\ep{\gamma_k}}$ between 1-forms.   It will be convenient to prove also that $F(\id_x)=1$ for any $x\in X$. Indeed,
\begin{equation*}
F(\id_x)= F_{\ep{\id_x}}(\tau) = \mathfrak{P}^{\infty}(A_{\ep{\id_x}})(\tau) =1  
\end{equation*}
since $A_{\ep{\id_x}}= \ep{\id_x}^{*}A = 0$ because $\ep{\id_x}$ has rank zero (Lemmata \ref{lem:thinprop} (b) and \ref{lem:thinflat}).

\textit{2.) Smoothness}. Let $c\maps U \to PX$ be a plot, i.e. the map
\begin{equation*}
\alxydim{@C=1.3cm}{[0,1] \times U \ar[r]^-{\id \times c} & [0,1] \times PX \ar[r]^-{\ev} & X}
\end{equation*}
is smooth. We extend this map  to a plot $\ep c\maps  \R \times U \to X$. Consider   the inclusion $\iota_u\maps \R \to \R \times U$ defined by $\iota_u(t) \df  (t,u)$, and the associated smooth map
\begin{equation*}
\Gamma\maps U \to \pt(\R \times U)\maps  u \mapsto \pt\iota_u(\tau)\text{.}
\end{equation*}
It follows that
\begin{equation}
\label{eq:Fc}
\alxydim{}{U \ar[r]^-{\Gamma} & \pt (\R \times U)  \ar[r]^-{F_{\ep c}} & G}
\end{equation} 
is smooth. We claim that \erf{eq:Fc} coincides with $F \circ c$; this shows that $F$ is smooth. Indeed, 
\begin{equation*}
F_{\ep c}(\Gamma(u)) = F_{\ep c}(\pt\iota_u(\tau)) = \iota_u^{*}F_{\ep c}(\tau) =  F_{\ep {c(u)}}(\tau) = F(c(u))\text{.}
\end{equation*}

\textit{3.) Thin homotopy invariance}. We have to show that $F(\gamma_1)=F(\gamma_2)$ whenever there is a thin homotopy $h\in PPX$ between  $\gamma_1$ and $\gamma_2$. The adjoint map $\exd{h}\maps [0,1]^2 \to X$ can be extended to a plot $\ep{\exd h}\maps  \R^2 \to X$ \cite[Section 2.3]{schreiber5}. One checks that
\begin{equation*}
((\exd h)^{*}F)(\gamma) = F_{\ep{P\exd h(\gamma)}}(\tau)=\mathfrak{P}^{\infty}(A_{{\ep{P\exd h(\gamma)}}})(\tau) = \mathfrak{P}^{\infty}(\ep\gamma^{*}A_{\ep{\exd h}})(\tau) = \mathfrak{P}^{\infty}(A_{\ep{\exd h}})(\gamma)\text{.}
\end{equation*}
Since $\exd h$ has rank one, the 1-form $A_{\ep{\exd h}} =(\ep{\exd h})^{*}A $ is flat by Lemma \ref{lem:thinflat}, and so is $(\exd h)^{*}F$. Then, the claim follows from Lemma \ref{lem:flathom} below.

By now we have defined a functor $\mathfrak{P}\maps  \mathcal{Z}_X^1(G) \to \fun XG$. It remains to show that it is inverse to the functor $\mathfrak{D}$. Suppose $A$ is an object in $\mathcal{Z}_X^1(G)$, and $F \df  \mathfrak{P}(A)$.
Over a plot $c\maps U \to X$, we have to check that $\mathfrak{D}^{\infty}(c^{*}F)=A_c$.
This is equivalent to $c^{*}F= \mathfrak{P}^{\infty}(A_c)$. Indeed, we find for $\gamma \in PU$
\begin{equation*}
(c^{*}F)(\gamma) = F(P c (\gamma)) = \mathfrak{P}^{\infty}(A_{\ep{P c (\gamma)}})(\tau) \stackrel{(*)}{=} \mathfrak{P}^{\infty}(\ep\gamma^{*}A_{c})(\tau) = \ep\gamma^{*}\mathfrak{P}^{\infty}(A_{c})(\tau) = \mathfrak{P}^{\infty}(A_c)(\gamma)
\end{equation*}
where $(*)$ comes from the commutative diagram
\begin{equation*}
\alxydim{@C=0.4cm}{\R \ar[rr]^-{\ep\gamma} \ar[dr]_{\ep{Pc(\gamma)}} && U_{c} \ar[dl]^{c}\\&X&}
\end{equation*}
and the last equality comes from the fact that $(\mathcal{P}\ep\gamma)(\tau)= \gamma$ as elements in $\pt X$.
Conversely, suppose an object $F$ in $\fun XG$ is given and $A \df  \mathfrak{D}(F)$. We have to check that $\mathfrak{P}(A)=F$. Indeed, for $\gamma\in PX$,
\begin{equation*}
\mathfrak{P}(A)(\gamma) = \mathfrak{P}^{\infty}(A_{\ep\gamma})(\tau) = \mathfrak{P}^{\infty}(\mathfrak{D}^{\infty}(\ep\gamma^{*}F))(\tau)  = (\ep\gamma^{*}F)(\tau) =F(\mathcal{P}\ep\gamma(\tau)) = F(\gamma)\text{.}
\end{equation*}
This completes the proof.
\end{proof}

We remark that the condition that the functors $\mathfrak{D}$ and $\mathfrak{P}$ restrict to the functors $\mathfrak{D}^{\infty}$ and $\mathfrak{P}^{\infty}$ over each plot, implies that the isomorphism of Theorem \ref{th:formsfunctors} commutes with the pullback functors \erf{eq:pullbacks} and respects flatness.
The supplementary lemma used in the proof above and   in Section \ref{sec:connections} is 

\begin{lemma}
\label{lem:flathom}
Let $X$ be a diffeological space. Suppose $F\maps  PX \to G$ is a smooth map satisfying
\begin{equation*}
F(\gamma' \circ \gamma) = F(\gamma') \circ F(\gamma)
\quand 
F(\id_x)=1
\end{equation*}
for all composable paths $\gamma',\gamma$, and all points $x\in X$. Suppose further that $h \in PPX$ is a homotopy between paths $\gamma_1$ and $\gamma_2$, such that $(\exd h)^{*}F$ is flat, where $\exd h\maps [0,1]^2 \to X$ is the adjoint of $h$. Then, $F(\gamma_1)=F(\gamma_2)$.
\end{lemma}

\begin{proof}
Recall that $\gamma_1,\gamma_2\maps  [0,1] \to X$ are locally constant in neighborhoods $U_1,U_2$ of $\left \lbrace 0,1 \right \rbrace$, respectively. Let $U \df  [0,1] \setminus (U_1 \cap U_2)$. Choose a smoothing function $\varphi$  such that $\varphi(t)=t$ for all $t\in U$. As a consequence, $\gamma_i=\gamma_i \circ \varphi$ for $i=1,2$. We can regard $\varphi$ as a path in $\R$, and use the smooth maps $\iota_s^h\maps \R \to \R^2\maps  t \mapsto (s,t)$ and $\iota_s^v\maps \R \to \R^2\maps t \mapsto (t,s)$ to construct paths in $\R^2$ like shown in Figure \ref{fig:paths}.
\begin{figure}
\begin{equation*}
\alxydim{@C=1.5cm}{(0,0) \ar[r]^{P\iota^h_0(\varphi)} \ar[d]_{P\iota^v_0(\varphi)} & (0,1) \ar[d]^{P\iota^v_1(\varphi)} \\ (1,0) \ar[r]_{P\iota_1^h(\varphi)} & (1,1)}
\end{equation*}
\caption{Paths around the boundary of the unit square.}
\label{fig:paths}
\end{figure}
Notice that $F(\gamma_i) = F_h(P\iota_0^h(\varphi))$ for $i=1,2$, where we have denoted $F_h \df  (\exd h)^{*}F$. For $(x,y) \df  \ev(\gamma_1)=\ev(\gamma_2)$, we calculate
\begin{equation*}
F(\gamma_1) = F(\id_y) \cdot F(\gamma_1)
=
F(P(\exd h \circ \iota_1^v)(\varphi) \pcomp P(\exd h \circ \iota_0^h)(\varphi)) 
= F_h(P\iota^v_1(\varphi) \pcomp P\iota^h_0(\varphi)) 
\end{equation*}
Since the paths $P\iota^v_1(\varphi) \pcomp P\iota^h_0(\varphi)$ and $P\iota_1^h(\varphi) \pcomp P\iota^v_0(\varphi)$ are obviously homotopic in $\R^2$, the latter result is equal to
\begin{equation*}
F_h(P\iota_1^h(\varphi) \pcomp P\iota^v_0(\varphi))
=F(P(\exd h \circ \iota_1^h)(\varphi) \pcomp P(\exd h \circ \iota^v_0)(\varphi))
=
F(\gamma_2) \cdot F(\id_x)
= F(\gamma_2)
\end{equation*}
Both lines together show the claim.
\end{proof}
\end{appendix}

\kobib{../../bibliothek/tex}

\end{document}